\documentclass[12pt,reqno]{amsart}

\usepackage{amssymb}        
\usepackage{latexsym}       
\usepackage{graphicx}       
\usepackage{multirow}       
\usepackage{xcolor}         

\usepackage[all]{xy}
\usepackage{caption,stmaryrd}
\usepackage{setspace}

\usepackage{mathrsfs}
\usepackage{verbatim}
\usepackage{booktabs}
\usepackage{tikz}
\usetikzlibrary{positioning}

\newtheorem{theorem}{Theorem}[section]
\newtheorem{lemma}[theorem]{Lemma}
\newtheorem{proposition}[theorem]{Proposition}
\newtheorem{corollary}[theorem]{Corollary}

\theoremstyle{definition}
\newtheorem{definition}[theorem]{Definition}
\newtheorem{problem}[theorem]{Problem}

\newtheorem{example}[theorem]{Example}

\newtheorem{hypothesis}[theorem]{Hypothesis}

\usepackage{comment}

\usepackage[backref=page,hypertexnames=false]{hyperref}

\usepackage{etoolbox}
\newcommand{\declaremathfont}[2]{%
    \begingroup
    \newcommand{\declareone}[1]{%
        \expandafter\gdef\csname #1##1\endcsname{#2{##1}}%
    }%
    \forcsvlist{\declareone}{A,B,C,D,E,F,G,H,I,J,K,L,M,N,O,P,Q,R,S,T,U,V,W,X,Y,Z}%
    \endgroup
}
\declaremathfont{bb}{\mathbb}
\declaremathfont{cal}{\mathcal}
\declaremathfont{rm}{\mathrm}
\declaremathfont{bf}{\mathbf}
\declaremathfont{sf}{\mathsf}
\declaremathfont{tt}{\mathtt}
\declaremathfont{it}{\mathit}
\declaremathfont{frak}{\mathfrak}

\newcommand{\declaremathgroup}[1]{\expandafter\newcommand\csname #1\endcsname{{\mathrm{#1}}}}
\forcsvlist{\declaremathgroup}{%
    SL,GL,Sp,SU,SO,CSp,
    PSL,PGL,PSp,PSU,
    ASL,AGL,
    Ree,Sz,Suz,Ru,Th,He,McL,HN,HS,Co,
}

\newcommand{\declaremathobject}[1]{\expandafter\newcommand\csname #1\endcsname{{\mathrm{#1}}}}
\forcsvlist{\declaremathobject}{%
    PG,
    AF,HA,AS,PA,TW,SD, CD, HC,
    Af,Ha,As,Pa,Tw,Sd,Cd,
}

\newcommand{\declaremathoperator}[1]{\expandafter\newcommand\csname #1\endcsname{{\mathsf{#1}}}}
\forcsvlist{\declaremathoperator}{%
    Ker,Tr,tr,rad,Rad,Soc,soc,Aut,Inn,Out,Hom,End,Syl,diag,Cay,Cos,Ann,rank,cl,core,Sym,
}

\def\a{{\alpha}}    {
}  \def\o{{\omega}} 
\def\Ga{{\Gamma}}

\def\ZZ{\mathbb{Z}}    

 \def\lcm{{\rm lcm}} 
  
   \def\geq{\geqslant}
\def\ov{\overline} 
\def\l{\langle} \def\r{\rangle}

 \def\Core{{\sf Core}} \def\Cos{{\sf Cos}}

\def\Aut{{\mathrm{Aut}}}  \def\Out{{\mathrm{Out}}}  
\def\GL{{\mathrm{GL}}}

\def\D{\mathrm{D}}    
\def\K{\mathrm{K}}    
\def\S{\mathrm{S}}    

\def\K{{\bf K}}

\def\RevMap{\mathsf{RevMap}}
\def\BiRevMap{\mathsf{BiRevMap}}
\def\Rev{\mathsf{Rev}}
\def\BiRev{\mathsf{BiRev}}

\begin{document}

\title[Coprime vertex-reversing maps]{A classification of vertex-reversing maps with Euler characteristic coprime to the edge number}

\author{Luyi Liu}
\author{Hanyue Yi}
\author{Yan Zhou Zhu}
\address{L. Liu, Shenzhen International Center for Mathematics \\
    Southern University of Science and Technology\\
    Shenzhen 518055\\
    P. R. China}
\email{12031108@mail.sustech.edu.cn}

\address{H. Yi, Department of Mathematics\\
    Southern University of Science and Technology \\
    Shenzhen 518055\\
    P. R. China}
\email{12431017@mail.sustech.edu.cn}

\address{Y.Z. Zhu, School of Mathematical Sciences \\
    Xiamen University \\
    Xiamen 361005\\
    P. R. China}
\email{zhuyz@xmu.edu.cn}

\date\today

\begin{abstract}
  An arc-regular map is \emph{vertex-reversing} if the automorphism group has dihedral vertex stabilizers.
    This paper classifies solvable vertex-reversing maps whose Euler characteristic is coprime to the edge number.
    The classification establishes that such maps fall into four families:
    $\D_{2n}$-maps, $(\D_{2m}\times\D_{2n})$-maps,
    $(\ZZ_{mn\ell}{:}\D_4)$-maps, and $(\ZZ_{3^f n}.\S_4)$-maps,
    with the parameters specified in the main theorem.
    Moreover, for each family, we provide an explicit formula for the Euler characteristic.

\end{abstract}
\maketitle

\section{Introduction}
A \textit{map} $\calM=(V,E,F)$ is a $2$-cell embedding of a graph $\Gamma$ in a closed surface $\calS$, where $V$, $E$, and $F$ are the sets of vertices, edges, and faces, respectively.
Throughout the paper, we assume that $|V|\geqslant 3$ and $|F|\geqslant 3$, and allow multiple edges but no loops or semi-edges.
The \textit{Euler characteristic} $\chi$ of $\calM$ is defined as $|V|-|E|+|F|$, which is equal to the Euler characteristic of the surface $\calS$.
An \textit{arc} of $\calM$ is an incident pair $(\a,e)$ consisting of a vertex $\a$ and an edge $e$, and a \textit{flag} is a mutually incident triple $(\a,e,f)$ consisting of a vertex $\a$, an edge $e$, and a face $f$.
An automorphism of $\calM$ is a permutation of the flag set that preserves incidence; the set of all automorphisms of $\calM$ forms the automorphism group, denoted by $\Aut(\calM)$.
The map $\calM$ is \textit{regular} if $\Aut(\calM)$ acts regularly on its flag set.
For a subgroup $G \leqslant \Aut(\calM)$, the map $\calM$ is called \textit{$G$-edge-transitive} (respectively, \textit{$G$-arc-transitive}) if $G$ acts transitively on the set of edges (respectively, arcs) of $\calM$.
It is called \textit{$G$-edge-regular} (respectively, \textit{$G$-arc-regular}) if $G$ acts regularly on the set of edges (respectively, arcs) of $\calM$.
We omit $G$ when the group $G$ is clear from the context.

Highly symmetric maps have been extensively studied from three perspectives: underlying graphs~\cite{LP2025}, supporting surfaces~\cite{ThmMOri}, and automorphism groups~\cite{Brahana,CFLZ,Coxeter,etm14tp}.
Regarding classifications based on Euler characteristic, numerous results exist for regular maps~\cite{regular-p,regular-3p,ori-reg--2p,Maps-sqfreeEC,ori-reg--2p2,reg--p3}, whereas results for edge-transitive maps remain rare.
The stabilizers of vertices, edges, and faces in the automorphism group of a map are cyclic or dihedral; see~\cite{Auto-ETmaps} for detailed proofs.
Using this property, Lemma~\ref{lem:Sylow} shows that, for a map $\calM=(V,E,F)$ with Euler characteristic $\chi$ and $\gcd(|E|,\chi)=1$, the automorphism group $\Aut(\calM)$ has cyclic or dihedral Sylow subgroups.
This naturally motivates one to study the following classification problem:
\begin{problem}\label{prob:1}
    Characterize and classify edge-transitive maps with Euler characteristic coprime to the edge number.
\end{problem}

The fourteen types of edge-transitive maps were first introduced by Graver and Watkins~\cite{etm14tp}.
Among them, one type is flag-regular, which has been extensively studied; four types are arc-regular, and these five types comprise the class of arc-transitive maps.
Furthermore, arc-regular maps are divided into two categories, where $\a$ is a vertex and $E(\a)$ is the set of edges incident with $\a$:
\begin{enumerate}
\item[(i)] {\it vertex-rotary} if $G_\alpha$ is cyclic and acts regularly on $E(\a)$;
\item[(ii)] {\it vertex-reversing} if $G_\alpha$ is dihedral and acts regularly on $E(\a)$.
\end{enumerate}
See~\cite{RevMap,RotaMap} for a theory of such maps and the associated construction methods.

To address Problem~\ref{prob:1}, we first focus on arc-transitive maps.
In subsequent work~\cite{Rotary}, we shall study vertex-rotary maps.
For vertex-reversing maps, the case with non-solvable automorphism groups was solved in~\cite{non-solvable}; hence the present paper focuses on the case with solvable automorphism groups and gives the classification in Theorem~\ref{thm:maps}.
We remark that Theorem~\ref{thm:maps} applies to maps with at least three vertices and three faces; see~\cite{fewV&F,ET-uniface,Li-Siran,uniface-dessins} for characterizations of maps with few vertices or faces.
	
To state the main result, we note that $G$-vertex-reversing maps comprise two distinct subclasses (see Subsection~\ref{subsec:defs-VerRevMap}): $G$-reversing maps, corresponding to Graver--Watkins type $2^*$, and $G$-bireversing maps, corresponding to Graver--Watkins type $2^{\rm P}$.
We denote by $\bfC_n^{(\lambda)}$ the multigraph obtained from a cycle on $n$ vertices by replacing each edge with $\lambda$ parallel edges. 
Similarly, $\K_n^{(\lambda)}$ denotes the multigraph obtained from the complete graph on $n$ vertices by replacing each edge with $\lambda$ parallel edges.
Thus $\bfC_n^{(\lambda)}$ has vertex-valency $2$ and edge-valency $2\lambda$.
When $\lambda=1$, $\bfC_n^{(1)}=\bfC_n$ and $\K_n^{(1)}=\K_n$.
For graphs $\Gamma_i=(V_i,E_i)$, $i=1,2$, their {\it direct product} $\Gamma_1\times\Gamma_2$ has vertex set $V_1\times V_2$, with $(u_1,u_2)$ adjacent to $(v_1,v_2)$ if and only if $\{u_i,v_i\}\in E_i$ for both $i=1,2$; and their {\it lexicographic product} $\Gamma_1[\Gamma_2]$ has vertex set $V_1\times V_2$, with $(u_1,u_2)$ adjacent to $(v_1,v_2)$ if and only if either $\{u_1,v_1\}\in E_1$, or $u_1=v_1$ and $\{u_2,v_2\}\in E_2$.
Throughout the paper, $\D_{2n}$ denotes the dihedral group of order $2n$.

\begin{theorem}\label{thm:maps}
    Assume that $\calM=(V, E, F)$ is a $G$-vertex-reversing map with solvable $G$ and $|V|,|F|\geqslant 3$.
    Let $\chi$ be the Euler characteristic of $\calM$ and let $\Gamma$ be the underlying graph.
    If $\gcd(\chi,|E|)=1$, then one of the following holds.
    \begin{enumerate}
        \item $\calM$ is a $\D_{2n}$-bireversing map with $\Gamma\cong\bfC_{d}^{(n/d)}$, and there exist integers $j$, $k$ such that $\gcd(j,k,n)=1$, $d=\gcd(j,n)\geqslant 3$ and 
        \[\chi=\gcd(j,n)+\gcd(2k-j,n)-n.\]
        In addition, $\calM$ is orientable if and only if $\gcd(j,2k,n)=2$.
        \item $\calM$ is a $G$-reversing map in Table~\ref{tab:classification}, where $m,n,\ell$ are pairwise coprime odd integers, $f\geqslant 0$ and $\gcd(6,s)=1$.
    \end{enumerate}
\end{theorem}
\begin{table}[htbp]
    \centering
    \renewcommand{\arraystretch}{1.2}
    \caption{$G$-reversing maps with $\gcd(\chi,|E|)=1$.}
    \label{tab:classification}
    \begin{tabular}{cccccc}
        \toprule
        $G$ & $\Gamma$& $\chi$ & {\rm orientable?}& {\rm ref}\\
        \midrule
        \multirow{2}{*}{$\D_{2n}$} & $\bfC_{n/2}^{(2)}$ & $1$ & $\times$&{\rm Thm~\ref{thm:type-I}\,(1)}\\
        \addlinespace[2pt]
        & $\bfC_{d}^{(n/d)}$ & {\rm see Thm~\ref{thm:type-I}\,(2)} & $\checkmark$&{\rm Thm~\ref{thm:type-I}\,(2)}\\
            \midrule
        \addlinespace[2pt]
        \multirow{2}{*}{$\D_{2m}\times\D_{2n}$}  & $\bfC_{m}^{(2n)}$& \multirow{2}{*}{$m+n-mn$}  &$\times$&\multirow{2}{*}{\rm Thm~\ref{thm:type-II}}\\
        \addlinespace[2pt]
        & $\bfC_{m}\times \bfC_n$&  &$\times$&\\
        \midrule
        \addlinespace[2pt]
        $\ZZ_{mn\ell}{:}\ZZ_2^2$ & $(\bfC_m\times \bfC_n)^{(\ell)}$ & $mn+m\ell+n\ell-2mn\ell$ &$\times$&{\rm Thm~\ref{thm:type-III}}\\
        \midrule
        \addlinespace[2pt]
        \multirow{3}{*}{$\ZZ_{3^{f}\cdot s}.\S_4$} & $\K_{4}^{(2\cdot3^{f}\cdot s)}$ & \multirow{3}{*}{$4-3^{f+1}\cdot s$} &$\times$&\multirow{3}{*}{\rm Thm~\ref{thm:type-IV}}\\
        \addlinespace[2pt]
        & $\bfC_{3^{f+1}s}[\overline{\K_2}]$ &  &$\times$&\\
        \addlinespace[2pt]
        & $\bfC_{3^{f+1}s}^{(4)}$ &  &$\times$&\\
        \bottomrule
    \end{tabular}
\end{table}

As shown in Lemma~\ref{lem:Sylow}, the coprime condition forces the automorphism group to have only cyclic or dihedral Sylow subgroups.
The following theorem determines all such solvable groups that can arise.
These groups are solvable and \textit{almost Sylow-cyclic} (each Sylow subgroup of odd order is cyclic and each Sylow $2$-subgroup contains a cyclic subgroup of index $2$).
Solvable almost Sylow-cyclic groups were classified in full generality by Zassenhaus~\cite{Zassenhaus}, but that classification is intricate and not well suited to our concrete setting.
Instead, in Section~3 we give a direct and self-contained proof of the following theorem that does not rely on Zassenhaus’ results; the four types are derived by elementary arguments from the structure of the maps.
A more explicit description of these groups in terms of generators will be given in Theorem~\ref{thm:Aut(M)}.

\begin{theorem}\label{thm:G}
    Assume that $\calM$ is a $G$-vertex-reversing map, where $G$ is a solvable group whose Sylow subgroups are cyclic or dihedral.   
    Then one of the following holds:
    \begin{enumerate}
        \item $G$ is a dihedral group;
        \item $G\cong (\ZZ_m\times\D_{2^e}){:}\ZZ_2\cong \D_{2^e}{:}\D_{2m}$, where $m$ is odd and $e\geqslant 2$;
        \item $G\cong \ZZ_{mn\ell}{:}\D_{2^e} \cong \ZZ_{2^{e-2}m}.(\D_{2n}\times \D_{2\ell})$, where $m,n,\ell$ are pairwise coprime odd integers and $e\geqslant 2$;
        \item $G\cong (\D_4\times\ZZ_n){:}\D_{2\cdot 3^{f+1}}\cong \ZZ_{3^fn}.\S_4$, where $\gcd(n,6)=1$ and $f\geqslant 0$.
    \end{enumerate}
\end{theorem}

It is worth noting that if the restriction on Sylow $2$-subgroups in Theorem~\ref{thm:G} is removed, then additional groups may occur; for instance, $\GL_2(3)\cong \rmQ_8{:}\S_3$ is an almost Sylow-cyclic whose Sylow $2$-subgroups are isomorphic to $\SD_{16}$, and it gives rise to vertex-reversing maps, as shown in Example~\ref{exam:Q_8.D_n-maps}.
However, by Lemma~\ref{lem:Sylow}, such groups do not appear in the main classification theorem, Theorem~\ref{thm:maps}.

The paper is organized as follows.
In Section~\ref{sec:Rev,BiRev}, we set up the necessary preliminaries on vertex-reversing maps and establish some basic properties that will be used throughout the paper.
The proof of Theorem~\ref{thm:G} is given in Section~\ref{sec:group}.
Finally, in Section~\ref{sec:final}, we analyze the four types obtained by reduction from Theorem~\ref{thm:4Tp-coprime} and complete the proof of Theorem~\ref{thm:maps}.


\section{Preliminaries on vertex-reversing maps}\label{sec:Rev,BiRev}

In this section, we recall the group-theoretic construction of vertex-reversing maps via arc-regular triples, and establish the key restriction that the coprime condition imposes on the automorphism group.

\subsection{Arc-regular triples and vertex-reversing maps}\label{subsec:defs-VerRevMap}\

Before presenting the group-theoretic construction of vertex-reversing maps, we describe their underlying graphs via coset graphs. 
The underlying graph of a vertex-reversing map is arc-transitive and may have multiple edges.
To handle such graphs, we use the generalized coset graph construction of~\cite[Construction~2.1]{RotaMap}, which allows arbitrary edge multiplicities.

\begin{definition}\label{def:graph}
    Let $H$ and $J$ be subgroups of the finite group $G$ such that $|J:H\cap J|=2$.
    Define an incidence structure $\Cos(G,H,J)=(V,E,\bfI)$ with
    \[
    V=[G:H],\quad E=[G:J]\quad\mbox{and}\quad \bfI=\{(Hx,Jy) : yx^{-1}\in JH\}.
    \]
    This is called a {\it coset graph} (allowing multiple edges).
  If $\Cos(G,H,J)$ is a simple graph, then it coincides with the classical coset graph $\Cos(G,H,HgH)$, where $g\in J{\setminus} H$, and is called the {\it base graph} of $\Cos(G,H,J)$.
\end{definition}

The following properties of coset graphs are taken from~\cite[Theorem~2.2]{RotaMap}.
Here, the vertex-valency is the number of adjacent vertices, and the edge-valency is the number of incident edges (counting multiplicities).

\begin{proposition}\label{prop:CosGraph}
    Let $\Ga=\Cos(G,H,J)$ with the notation of Definition~\ref{def:graph}.
    Take $g\in J\setminus (H\cap J)$, set $K=H\cap H^g$ (so that $H\cap J\leqslant K\leqslant H$), $L=\l K,g\r$, $k=|H:K|$, and $\lambda=|K:H\cap J|$.
    Suppose that $\l H,J\r=G$. 
    Then $\Ga$ is a connected $G$-arc-transitive graph and the following statements hold.
    \begin{enumerate}
        \item $\Ga$ has constant vertex-valency $k$, edge-multiplicity $\lambda$, and edge-valency $k\lambda$;
        \item $H$ is a vertex-stabilizer, $J$ is an edge-stabilizer, and $H\cap J$ is an arc-stabilizer;
        \item $\Ga$ is simple if and only if $L=J$.
    \end{enumerate}
\end{proposition}

Having described the underlying graphs via coset graphs, we now define vertex-reversing maps. Recall that we assume $|V|,|F|\geqslant 3$ throughout.

\begin{definition}\label{def:RevMap}
    An ordered triple of involutions $(x,y,z)$ of a group $G$ is called an {\it arc-regular triple} if $G=\l x,y,z\r$ and $|G:\l x,y\r|\geqslant 3$.
    Given such a triple, the coset graph $\Ga=\Cos(G,\l x,y\r,\l z\r)$ admits the following two embeddings:
    \begin{enumerate}
        \item $\RevMap(G,x,y,z)$, a reversing map with face set $F=[G:\l x,z\r]\cup [G:\l y,z\r]$;
        \item $\BiRevMap(G,x,y,z)$, a bireversing map with face set $F=[G:\l x,y^z\r]$.
    \end{enumerate}
    In both maps, a vertex $\a\in V$ or an edge $e\in E$ is incident with a face $f\in F$ if $\a\cap f\neq\emptyset$ or $e\cap f\neq\emptyset$, respectively.
\end{definition}
Note that $|G:\l x,y\r|\geqslant 3$ ensures $\l x,z\r\neq\l y,z\r$; otherwise $|G:\l x,y\r|$ would be $1$ or $2$. 
Hence, the two orbits in the face set of $\RevMap(G,x,y,z)$ are always distinct, so the definition is unambiguous.

By~\cite{RevMap}, every $G$-vertex-reversing map with at least three vertices and three faces is isomorphic to $\RevMap(G,x,y,z)$ or $\BiRevMap(G,x,y,z)$ for some arc-regular triple $(x,y,z)$ for $G$.
We say two triples are \textit{equivalent} if they admit isomorphic (bi-)reversing maps.
The following proposition records when two triples are equivalent.
\begin{proposition}\label{prop:iso-VerRevMaps}{\rm\cite{RevMap}}
    Two arc-regular triples $(x,y,z)$ and $(x',y',z')$ for $G$ give isomorphic (bi-)reversing maps if and only if $(x,y,z)^\sigma=(x',y',z')$ or $(y',x',z')$ for some $\sigma\in\Aut(G)$.
\end{proposition}

For a $G$-vertex-reversing map, the subgroup $G^+$ defined in~\cite[Table~2]{Auto-ETmaps} consists precisely of the orientation-preserving automorphisms when the map is orientable. In this case, $G^+$ has index $2$ in $G$~\cite{Auto-ETmaps}, yielding the following criterion.

\begin{lemma}\label{cri-Orient}
    A $G$-vertex-reversing map $\calM$ is orientable if and only if $|G:G^+|=2$, where
    $G^+=\begin{cases}
        \l xy,xz\r, & \text{if }\calM\cong \RevMap(G,x,y,z), \\
        \l xy,z,xzy\r, & \text{if }\calM\cong \BiRevMap(G,x,y,z).
    \end{cases}$
\end{lemma}

For later use, we record a lemma that characterises when the coset graph decomposes as a direct product of two cycles.

\begin{lemma}\label{lem:product-coset-graph}
    For $i=1,2$, let $G_i=\l a_i\r{:}\l u_i\r\cong\D_{2m_i}$, $H_i=\l u_i\r$, and let $g_i=a_i u_i$, where
    $m_i\geq 3$.
    Then
    \[
        \Cos(G,H,\l g\r)\cong \Cos(G,H,HgH)\cong \bfC_{m_1}\times \bfC_{m_2},
    \]
    where $G=G_1\times G_2$, $H=H_1\times H_2=\{(h_1,h_2)\mid h_1\in H_1,\, h_2\in H_2\}$, and $g=(g_1,g_2)$.
\end{lemma}
\begin{proof}
    Let $\Gamma=\Cos(G,H,\l g\r)$ and $\Delta=\bfC_{m_1}\times\bfC_{m_2}$.
    The vertices of $\Delta$ are the pairs $(j_1,j_2)$ with $j_i\in\ZZ_{m_i}$, where
    $(j_1,j_2)$ is adjacent to $(j_1',j_2')$ if and only if
    $j_i'-j_i\equiv\pm1\pmod{m_i}$ for $i=1,2$.
    By Proposition~\ref{prop:CosGraph}, the graph $\Gamma$ is simple and coincides
    with $\Cos(G,H,HgH)$.
    Note that the vertices of $\Gamma$ can be uniquely written as $H(a_1^{j_1},a_2^{j_2})$ for $j_1\in\ZZ_{m_1}$ and $j_2\in\ZZ_{m_2}$.
    Define
    \[
        \Phi:H(a_1^{j_1},a_2^{j_2})\mapsto (j_1,j_2).
    \]
    Clearly, $\Phi$ is a bijection between the vertex sets of $\Gamma$ and $\Delta$.
    We next check that $\Phi$ preserves adjacency; then $\Phi$ is an isomorphism between $\Gamma$ and $\Delta$.

    Note that $H_i g_i H_i=H_i a_iH_i=H_i a_i\sqcup H_i a_i^{-1}$ for $i=1,2$.
    Then 
    \[\begin{aligned}
        HgH&= (H_1,H_2) (g_1,g_2)(H_1,H_2)=(H_1,H_2)\{a_1^{\epsilon_1},a_2^{\epsilon_2}\mid \epsilon_1,\epsilon_2=\pm1\}(H_1,H_2)\\
        &=\bigcup_{\epsilon_1,\epsilon_2=\pm1}H(a_1^{\epsilon_1},a_2^{\epsilon_2}).
    \end{aligned}\]
    Hence two vertices $H(a_1^{j_1},a_2^{j_2})$ and
    $H(a_1^{j_1'},a_2^{j_2'})$ are adjacent in $\Gamma$ if and only if
    $j_i'-j_i\equiv\pm1\pmod {m_i}$ for $i=1,2$.
    This is exactly the adjacency condition for their images under $\Phi$ in
    $\Delta=\bfC_{m_1}\times\bfC_{m_2}$.
    Therefore $\Phi$ is an isomorphism between $\Gamma$ and $\Delta$.
\end{proof}

\subsection{The Euler characteristic and the coprime condition}\label{subsec:|E|}\

Recall that the Euler characteristic of a map $\calM=(V,E,F)$ is $\chi(\calM)=|V|-|E|+|F|$.
Throughout the sequel, for a group $G$, we denote by $\pi(G)$ the set of prime divisors of $|G|$; for $p\in\pi(G)$, we denote by $G_p$ a Sylow $p$-subgroup of $G$.
The next lemma characterizes the Sylow subgroups of the automorphism group of a map with Euler characteristic coprime to the number of edges.

\begin{lemma}\label{lem:Sylow}
	Let $\calM=(V,E,F)$, $G\leqslant\Aut(\calM)$, and let $\chi=\chi(\calM)$.
	Assume that $|V|,|F|\geqslant 3$ and $\gcd(\chi,|E|)=1$.
    Then the following statements hold.
	\begin{enumerate}
		\item $\gcd(\chi,|G|)\in\{1,2,4\}$, and if $\gcd(\chi,|G|)\neq 1$, then $|E|$ is odd;
		\item each Sylow subgroup of $G$ is contained in a vertex, edge, or face stabilizer;
		\item each Sylow subgroup of $G$ is cyclic or dihedral;
		\item $|G|=\lcm\{|G_\o|: \o\in V\cup E\cup F\}$.
	\end{enumerate}
\end{lemma}
\begin{proof}
    It is known that $G$ acts semiregularly on the flag set $\Phi$ of $\calM$.
    Since each edge is incident with $4$ flags,, it follows that $|G|$ divides $|\Phi|=4|E|$.
    Thus, $\gcd(\chi,|G|)$ divides $\gcd(\chi,4|E|)=\gcd(\chi, 4)$ as $\gcd(\chi,|E|)=1$.
    This implies that $\gcd(\chi,|G|)\in\{1,2,4\}$.
    If $\gcd(\chi,|G|)\neq 1$, then $\chi$ is even, and hence $|E|$ is odd as $\gcd(\chi,|E|)=1$.
    This proves part~(1).

    Decompose $V$, $E$ and $F$ into $G$-orbits by
    \[V=V_1\cup\cdots\cup V_r,\quad E=E_1\cup\cdots\cup E_s,\quad F=F_1\cup\cdots\cup F_t.\]
    Let $v_i\in V_i$, $e_j\in E_j$ and $f_k\in F_k$.
    Then we have that 
    \[\chi=|V|-|E|+|F|=\sum_{i=1}^r |G:G_{v_i}|-\sum_{j=1}^s |G:G_{e_j}|+\sum_{k=1}^t |G:G_{f_k}|.\]
    Suppose that a Sylow $p$-subgroup $P$ of $G$ has no fixed points on $V\cup E\cup F$.
    Then, by conjugacy of Sylow $p$-subgroups, no stabilizer $G_{v_i}$, $G_{e_j}$, or $G_{f_k}$ contains a Sylow $p$-subgroup of $G$.
    Consequently, $p$ divides every term in the above expression for $\chi$, giving $p\mid\chi$.
    Hence, $p\mid\gcd(\chi,|G|)$, so $p=2$.
    This implies that $|E|$ is odd as $\gcd(\chi,|G|)\neq 1$ by part~(1).
    Then some $|G:G_{e_j}|$ is odd, so $G_{e_j}$ contains a Sylow $2$-subgroup of $G$. 
    Thus, $P$ fixes an edge, a contradiction.
    Hence, part~(2) holds.
    Part~(3) follows, since each stabilizer $G_\o$ is cyclic or dihedral by~\cite{Auto-ETmaps}.

    Clearly, $\lcm\{|G_\o| : \o\in V\cup E\cup F\}$ divides $|G|$, since each $|G_\o|$ divides $|G|$.
    Note that $|G|=\prod_{p\in\pi(G)}|G_p|$, and each $|G_p|$ divides $|G_\o|$ for some $\o\in V\cup E\cup F$ by part~(2).
    Thus, $|G|$ divides $\lcm\{|G_\o| : \o\in V\cup E\cup F\}$.
    Therefore, $|G|=\lcm\{|G_\o| : \o\in V\cup E\cup F\}$, proving part~(4).
\end{proof}

Lemma~\ref{lem:Sylow}\,(3) forces $G$ to be almost Sylow-cyclic and to have cyclic or dihedral Sylow $2$-subgroups.
Without the coprime condition, other types of almost Sylow-cyclic groups can also give rise to vertex-reversing maps, as illustrated below.

\begin{example}\label{exam:Q_8.D_n-maps}
Let $G=\GL_2(3)\cong \rmQ_8{:}\S_{3}$, and let $x=\begin{pmatrix} 1&0\\0&-1\end{pmatrix}$ and $r=\begin{pmatrix} -1&-1\\1&0\end{pmatrix}$.
Then $|x|=2$, $|r|=3$, and $G=\langle x,r\rangle=\langle x,x^r,x^{r^2}\rangle$.
Thus, $G$ is almost Sylow-cyclic, and its Sylow $2$-subgroups are isomorphic to $\SD_{16}$.
Note that $(x,x^r,x^{r^2})$ is an arc-regular triple for $G$ and
\begin{enumerate}
    \item $\RevMap(G,x,x^r,x^{r^2})$ has $4$ vertices, $24$ edges and $8$ faces with Euler characteristic $-12$;
    \item $\BiRevMap(G,x,x^r,x^{r^2})$ has $4$ vertices, $24$ edges and $6$ faces with Euler characteristic $-14$.
\end{enumerate}
\end{example}

For a $G$-vertex-reversing map $\calM$, let $(x,y,z)$ be an arc-regular triple for $\calM$. 
By Definition~\ref{def:RevMap}, we have
\[\chi(\calM)=\left\{\begin{aligned}
    &\frac{|G|}{|\l x,y\r|}-\frac{|G|}{2}+\frac{|G|}{|\l x,z\r|}+\frac{|G|}{|\l y,z\r|},&\mbox{ if $\calM\cong\RevMap(G,x,y,z)$;}\\
    &\frac{|G|}{|\l x,y\r|}-\frac{|G|}{2}+\frac{|G|}{|\l x,y^z\r|},&\mbox{ if $\calM\cong \BiRevMap(G,x,y,z)$.}
\end{aligned}\right.\]
Note that the above formula for a reversing map is independent of the ordering of the involutions $x,y,z$.
Thus, we immediately obtain the following useful lemma.
\begin{lemma}\label{ChiRev=ChiRev}
    Let $(x,y,z)$ be an arc-regular triple for $G$.
    If $\{x',y',z'\}=\{x,y,z\}$, then $\RevMap(G,x,y,z)$ and $\RevMap(G,x',y',z')$ have the same Euler characteristic.
\end{lemma}

\section{Automorphism groups}\label{sec:group}

This section classifies the solvable groups $G$ that can act vertex-reversingly on a map $\calM$ with $\gcd(\chi(\calM),|E|)=1$.
By Lemma~\ref{lem:Sylow}, if $G$ acts vertex-reversingly on a map with $\gcd(\chi,|E|)=1$,
then every Sylow subgroup of $G$ is cyclic or dihedral.
Thus, we always assume the following hypothesis holds in this section.
\begin{hypothesis}\label{hypo:Aut(M)}
    Let $G$ be a solvable group with an arc-regular triple $(x,y,z)$ such that each of its Sylow subgroups is cyclic or dihedral.
\end{hypothesis}

We will proceed in two steps:
\begin{enumerate}
    \item classify all solvable groups with this property (Theorem~\ref{thm:Aut(M)});
    \item apply the coprime condition to obtain the final four types (Theorem~\ref{thm:4Tp-coprime}).
\end{enumerate}

For the rest of this section, we always assume that $G$ and $(x,y,z)$ satisfy Hypothesis~\ref{hypo:Aut(M)}.
Since $G$ is generated by three involutions, each abelian quotient of $G$ is generated by involutions, and hence is elementary abelian.
\begin{lemma}\label{lem:G/N}
Let $N$ be a normal subgroup of $G$.
If $G/N\neq 1$ is abelian, then $G/N$ is an elementary abelian $2$-group.
\end{lemma}

If the Sylow $2$-subgroup of $G$ is also cyclic, that is, $G$ is \textit{Sylow-cyclic}, $G$ is known to be a metacyclic group with $G'$ and $G/G'$ as two cyclic groups of coprime orders, see~\cite[page 447]{Huppert2025} for instance.
In this case, we deduce that $G$ is dihedral.
\begin{corollary}\label{coro:Syl-cyc->dih}
    If $G$ is Sylow-cyclic, then $G$ is a dihedral group.
\end{corollary}
\begin{proof}
    In this case, $G=A{:}B\cong \bbZ_m{:}\bbZ_n$, where $A=G'\cong\ZZ_m$ and $B\cong\ZZ_n$ with $\gcd(n,m)=1$.
    Then $n=2$ and $m$ is odd by Lemma~\ref{lem:G/N}.
    Assume that $B=\l \beta\r$.
    Then $A=\l \alpha_0\r\times\l\alpha_1\r$, where $\gcd(|\alpha_0|,|\alpha_1|)=1$ and $(\alpha_0,\alpha_1)^\beta=(\alpha_0,\alpha_1^{-1})$.
    Hence, $G/\l\alpha_1\r\cong \l\alpha_0\r\times \l \beta\r$, and then it is an elementary abelian $2$-group by Lemma~\ref{lem:G/N}.
    It follows that $\alpha_0=1$, and so $G\cong \D_{2m}$.
\end{proof}

Recall that the Fitting subgroup $F$ of $G$ is the maximal nilpotent normal subgroup of $G$, and $\bfC_G(F)\leqslant F$ as $G$ is solvable.
Note that the Hall $2'$-subgroup $F_{2'}$ is cyclic as $F$ is also almost Sylow-cyclic.
Then we have the following reduction.
For a solvable group $G$ and a set of primes $\{p_1,\ldots,p_i\}$, we denote by $G_{\{p_1,\ldots,p_i\}}$ and $G_{\{p_1,\ldots,p_i\}'}$ some Hall $\{p_1,\ldots,p_i\}$-subgroup and Hall $\{p_1,\ldots,p_i\}'$-subgroup of $G$, respectively; these subgroups exist by Hall's theorem.
\begin{lemma}\label{lem:fit}
    Let $F$ be the Fitting subgroup of $G$.
    Then $F_{2'}$ is cyclic, and one of the following cases holds.
    \vspace{-.2cm} 
    \begin{spacing}{1.2} 
    \begin{enumerate}
        \item[\rm(A)] $G=F$ is a dihedral $2$-group;
        \item[\rm(B)] $F_{2'}=G_{2'}$ and $G/\bfC_{G}(G_{2'})\cong G_2/\bfC_{G_2}(G_{2'})\cong \ZZ_2$ or $\ZZ_2^2$;
        \item[\rm(C)] $F_{\{2,3\}'}=G_{\{2,3\}'}$, and $G_{\{2,3\}}\cong\ZZ_2^2{:}\D_{2\times 3^e}\cong (\ZZ_2^2\times \ZZ_{3^{e-1}}).\S_3
        \cong \ZZ_{3^{e-1}}.\S_4$ satisfying
        \[
        G/\bfC_{G}(G_{\{2,3\}'})\cong G_{\{2,3\}}/\bfC_{G_{\{2,3\}}}(G_{\{2,3\}'})\cong  \ZZ_2\mbox{ or }1.
        \]
    \end{enumerate}
    \end{spacing}
\end{lemma}
\vspace{-1.1cm} 
\begin{proof}
    If $G=F$, then $G$ is nilpotent and hence is a direct product of its Sylow subgroups.
    By Lemma~\ref{lem:G/N}, $G/G_2$ is trivial, and thus $G$ is a cyclic or dihedral $2$-group.
    Note that no cyclic group possesses arc-regular triples, and thus $G$ is a dihedral $2$-group as in Case~A.
    From now on, we assume that $F\neq G$.

    First, we assume that $G/\bfC_G(F_2)$ is a $2$-group; namely each odd order element acts trivially on $F_2$.
    Note that $G/F\lesssim\Out(F)$ as $G$ is solvable.
    Then there is a natural map $\varphi: G\rightarrow \Out(F)$ with kernel $F$.
    As $F$ is nilpotent, we have the following decomposition
    \[\Out(F)=E\times O\cong \Out(F_2)\times\Out(F_{2'}),\]
    where $E\cong\Out(F_2)$ acts trivially on $F_{2'}$ and $O\cong\Out(F_{2'})$ acts trivially on $F_{2}$.
    Note that $O\cong \Out(F_{2'})$ is abelian as $F_{2'}$ is cyclic.
    Then 
    \[\varphi(G)\leqslant \Out(F) =E\times O_2\times O_{2'}.\]
    Since $\varphi(G)=G/F$, it follows that the projection image of $\varphi(G)$ on $O_{2'}$ is also isomorphic to a quotient of $G$.
    Lemma~\ref{lem:G/N} shows that $G$ has no abelian quotient of odd order, and thus $G$ projects trivially on $O_{2'}$.
    Recall that $G/\bfC_G(F_2)$ is a $2$-group.
    It follows that $\varphi(G)_{2'}$ acts trivially on $F_2$, and so $\varphi(G)_{2'}\leqslant O_{2'}$.
    Hence we deduce that $\varphi(G)_{2'}=1$, that is, $G/F$ is a $2$-group.
    Thus $G_{2'}\leqslant F$, and so $G_{2'}=F_{2'}$ is cyclic.
    We have that $G_2/\bfC_{G_2}(G_{2'})\cong G/\bfC_{G}(G_{2'})$ is an abelian $2$-group as $G=G_{2'}{:}G_2$.
    If $G= G_{2'}\times G_2$, then $G_{2'}$ is trivial and $G=F$ by Lemma~\ref{lem:G/N}, leading to a contradiction.
    Thus $G_2/\bfC_{G_2}(G_{2'})\neq 1$ is an elementary abelian $2$-group by Lemma~\ref{lem:G/N}.
    Since $G_2$ is either cyclic or dihedral, we have that $G_2/\bfC_{G_2}(G_{2'})\cong \ZZ_2$ or $\ZZ_2^2$ as in Case~B.

    Now, we assume that $G/\bfC_G(F_2)$ is not a $2$-group.
    Then $F_2\cong\ZZ_2^2$ and $G/\bfC_G(F_2)\cong\S_3$ as $F_2$ is either cyclic or dihedral.
    This implies that $G_{\{2,3\}'}\leqslant\bfC_G(F_2)$.
    Since $F_{2'}$ is cyclic, $G/\bfC_G(F_{2'})$ is abelian, and then it is either trivial or an elementary abelian $2$-group by Lemma~\ref{lem:G/N}.
    Then $G_{2'}\leqslant \bfC_G(F_{2'})$; in particular $G_{\{2,3\}'}\leqslant \bfC_{G}(F_{2'})$.
    It follows from $\bfC_{G}(F)=\bfC_{G}(F_2)\cap \bfC_{G}(F_{2'})$ that $G_{\{2,3\}'}\leqslant \bfC_G(F)\leqslant F$, and so $G_{\{2,3\}'}=F_{\{2,3\}'}\lhd G$.
    Note that $G_{\{2,3\}}/F_2$ is Sylow-cyclic.
    Then $G_{\{2,3\}}/F_2\cong\D_{2\times 3^e}$.
    Since $G/\bfC_{G}(F_2)\cong \S_3$ is not abelian, we have that 
    \[G_{\{2,3\}}/F_2\cong \D_{2\times 3^e}\mbox{ and }G_{\{2,3\}}/\bfC_{G}(F_2)\cong \S_3.\]
    This implies that 
    \[G_{\{2,3\}}\cong \ZZ_2^2.\D_{2\times 3^e}\cong (\ZZ_2^2\times \ZZ_{3^{e-1}}).\S_3\cong \ZZ_{3^{e-1}}.\S_4.\]
    We now show that the first extension splits.
    Let $K=(G_{\{2,3\}})'$.
    Then 
    \[K=F_2{:}K_3\cong \ZZ_2^2{:}\ZZ_{3^e}\cong (\ZZ_2^2\times \ZZ_{3^{e-1}}).\ZZ_3\mbox{ has index $2$ in $G_{\{2,3\}}$}.\]
    Since $G_2\cong \D_8$, we can take an $x\in G_2\setminus F_2$.
    Then $K_3^x$ is also a Sylow $3$-subgroup of $K$, and $K_3^x=K_3^\alpha$ for some $\alpha\in F_2$.
    Let $x'=x\alpha^{-1}$.
    Then $x'\in G_2\setminus F_2$ and $x'$ normalizes $K_3$.
    Note that $(x')^2$ centralizes $K_3$ and $\bfC_{K_3}(F_2)=1$.
    Then $x'$ has order $2$, and hence
    \[G_{\{2,3\}}=F_2{:}(K_3{:} \langle x'\rangle)\cong\ZZ_2^2{:}\D_{2\times 3^e}.\]
    Finally, we show that $G/\bfC_{G}(G_{\{2,3\}'})\cong G_{\{2,3\}}/\bfC_{G_{\{2,3\}}}(G_{\{2,3\}'})\cong  \ZZ_2$.
    Since $G_{\{2,3\}'}=F_{\{2,3\}'}$ is cyclic, it follows that $G/\bfC_{G}(G_{\{2,3\}'})\cong G_{\{2,3\}}/\bfC_{G_{\{2,3\}}}(G_{\{2,3\}'})$ is abelian.
    Recall that $K$ is the commutator subgroup of $G_{\{2,3\}}$ and $|G_{\{2,3\}}:K|=2$.
    Thus $G/\bfC_{G}(G_{\{2,3\}'})\cong \ZZ_2$ or $1$ as in Case~C.
\end{proof}

The following is one of the main results of this section, which refines Theorem~\ref{thm:G} by giving an explicit description of the groups satisfying Hypothesis~\ref{hypo:Aut(M)}.
\begin{theorem}\label{thm:Aut(M)}
    Assume that $G$ is a solvable group satisfying Hypothesis~\ref{hypo:Aut(M)}.
    Then $G_{2'}=\l g\r$ is cyclic, $G_2=\l u\r{:}\l v\r$ with $|v|=2$, and one of the following holds.
    \begin{enumerate} 
        \item $G$ is a dihedral group;
        \item $G=\l g\r{:}\l u,v\r=(\l g\r\times \l u^2,v\r){:}\l uv\r$, where $g^u=g^{-1}$;
        \item $G=\l g\r{:}\l u,v\r=(\l g\r\times \l u^2\r).\l \ov u, \ov v\r$, and there is a decomposition $\l g\r=\l a\r\times\l b\r\times\l c\r$ with $bc\not=1$ such that $(a,b,c)^u=(a^{-1},b^{-1},c)$, and $(a,b,c)^v=(a^{-1},b,c^{-1})$.
        \item $G=(\l u^2,uv\r\times\l g_{3'}\r){:}\l g_3,v\r$,
        where $(g_{3'}g_3)^v=(g_{3'}g_3)^{-1}$, and $G/\l g_{3'}g_3^3\r\cong\S_4$.
    \end{enumerate}
\end{theorem}
\begin{proof}
    If $G$ is Sylow-cyclic, then $G$ is dihedral as in part~(1).
    We assume from now on that $G$ is not dihedral and has dihedral Sylow $2$-subgroups.
    Then one of Cases~B and C in Lemma~\ref{lem:fit} holds.
    Let $F$ be the Fitting subgroup of $G$.

    \textbf{Case B}: $F_{2'}=G_{2'}\lhd G$ and $G/\bfC_G(G_{2'})\cong G_2/\bfC_{G_2}(G_{2'})\cong\ZZ_2$ or $\ZZ_2^2$.

    Note that $G_{2'}=F_{2'}$ is cyclic as it is nilpotent with cyclic Sylow subgroups.
    Let $G_{2'}=\langle g\rangle$, and let $G_2=\l u\r{:}\l v\r\cong \D_{2^{f+1}}$.
    Then 
    \[G=G_{2'}{:}G_{2}=\langle g\rangle{:}\l u, v\r\cong \langle g\rangle{:}\D_{2^{f+1}}.\]
    Note that $\bfC_{G_2}(G_{2'})\geqslant (G_2)'=\langle u^2\rangle$.
    This implies that
    \[\bfC_{G_2}(G_{2'})=\l u^2\r,\ \l u\r,\ \l u^2\r{:}\l v\r\mbox{ or }\l u^2 \r{:}\l uv\r.\]
    If $\bfC_{G_2}(G_{2'})=\l u^2 \r{:}\l uv\r$, then we have $\bfC_{G_2}(G_{2'})=\l u^2\r{:}\l v_0\r$ and $G_2=\l u\r {:}\l v_0\r$ by setting $v_0=uv$.
    Hence it suffices to consider the cases where $\bfC_{G_2}(G_{2'})=\l u^2\r$, $\l u\r$ or $\l u^2\r{:}\l v\r$.

    Assume that $\bfC_{G_2}(G_{2'})=\l u\r$.
    Then $G_2=\bfC_{G_2}(G_{2'}){:}\l v\r $, and hence 
    \[G=G_{2'}{:}G_{2}=(G_{2'}\times \bfC_{G_2}(G_{2'})){:}\langle v\rangle=(\langle g\rangle\times \l u\r){:}\langle v\rangle.\]
    Note that $\langle u\rangle\lhd G$ and $G/\langle u\rangle\cong \langle g\rangle{:}\langle v\rangle$ is Sylow-cyclic.
    Then $ \langle g\rangle{:}\langle v\rangle$ is dihedral, so is $G=(\langle g\rangle\times \l u\r){:}\langle v\rangle=\langle gu\rangle{:}\langle v\rangle$, a contradiction.

    Assume that $\bfC_{G_2}(G_{2'})=\l u^2\r{:}\l v\r$.
    Then $G_2=\bfC_{G_2}(G_{2'}){:}\l uv\r $, and hence 
    \[G=( G_{2'}\times \bfC_{G_2}(G_{2'})){:}\l uv\r =(\l g\r \times \l u^2,v\r){:}\l uv\r.\]
    By arguments similar to those in the previous case, $G/\bfC_{G_2}(G_{2'})\cong \l g\r{:}\l uv\r$ is dihedral, and so $g^{uv}=g^{-1}$.
    Since $[g,v]=1$, we have that $g^{u}=g^{-1}$ as in part~(2).

    Assume that $\bfC_{G_2}(G_{2'})=\l u^2\r$.
    Then 
    \[G=(G_{2'}\times \bfC_{G_2}(G_{2'})).\l \ov u,\ov v\r =(\l g\r \times \l u^2\r).\l \ov u,\ov v\r.\]
    Note that each of $u$, $v$ and $uv$ admits an automorphism of $\l g\r$ of order at most $2$.
    Then $\l g\r$ has the following three coprime decompositions:
    \[
    \l g\r=\l\alpha_1\r\times \l\beta_1\r=\l\alpha_2\r\times \l\beta_2\r=\l\alpha_3\r\times \l\beta_3\r ,
    \]
    where $(\alpha_1,\beta_1)^{u}=(\alpha_1^{-1},\beta_1)$, $(\alpha_2,\beta_2)^{v}=(\alpha_2^{-1},\beta_2)$ and $(\alpha_3,\beta_3)^{uv}=(\alpha_3^{-1},\beta_3)$.
    Let
    \[
    \pi_a=\pi(\l \alpha_1\r)\cap \pi(\l \alpha_2\r),\ \pi_b=\pi(\l \alpha_1\r)\cap \pi(\l \alpha_3\r),\ \pi_c=\pi(\l \alpha_2\r)\cap \pi(\l \alpha_3\r).
    \]
    Note that $\pi_a\cap \pi_b=\pi(\l\alpha_1\r)\cap \pi(\l\alpha_2\r)\cap \pi(\l \alpha_3\r)$.
    Then
    \[(g_p)^{u}=(g_p)^{v}=(g_p)^{uv}=g^{-1}_p\mbox{ for every $p\in\pi_a\cap \pi_b$}.\]
    It follows that $g_p=1$ as $p\in\pi(\l g\r)$ is odd, and so $\pi_a\cap \pi_b=\emptyset$.
    Similarly, we have that $\pi_a\cap \pi_c=\pi_b\cap\pi_c=\emptyset$.
    Thus $\pi_a\cup \pi_b\cup \pi_c$ is a disjoint union.
    If this union does not equal $\pi(\l g\r)$, then
    \[ [g_p,\l u,v\r]=1\mbox{ for any $p\in\pi(\l g\r)\setminus (\pi_a\cup \pi_b\cup \pi_c)$}.\]
    This implies that $\bfC_{\l g\r}(G_2)\neq 1$.
    Since $\l g\r=\l [g, G_2]\r\times \bfC_{\l g\r}(G_2)$, we have that 
    \[G= \l g\r{:}\l  u, v\r=\bfC_{\l g\r}(G_2)\times \l [g, G_2]\r{:}\l u,v\r.\]
    By Lemma~\ref{lem:G/N}, $\bfC_{\l g\r}(G_2)\leqslant \l g\r$ is an elementary abelian $2$-group, which is impossible.
    Thus  $\pi(\l g\r)$ equals the disjoint union $\pi_a\cup \pi_b\cup \pi_c$.
    Then  $\l g\r$ has a decomposition
    \[\l g\r=\l g\r_{\pi_a}\times \l g\r_{\pi_b}\times \l g\r_{\pi_c}=\l a\r\times \l b\r\times \l c\r,\mbox{ where}\] 
    $(a,b,c)^u=(a^{-1},b^{-1},c)$ and $(a,b,c)^v=(a^{-1},b,c^{-1})$.
    Moreover, if $bc=1$ then $\bfC_{G}(\langle g\rangle)$ contains $uv$, a contradiction.
    Therefore, part~(3) holds.

    \textbf{Case C}: $F_{\{2,3\}'}=G_{\{2,3\}'}$, and $G_{\{2,3\}}\cong\ZZ_2^2{:}\D_{2\times 3^e}\cong (\ZZ_2^2\times \ZZ_{3^{e-1}}).\S_3$ satisfying
    \[
    G/\bfC_{G}(G_{\{2,3\}'})\cong G_{\{2,3\}}/\bfC_{G_{\{2,3\}}}(G_{\{2,3\}'})\cong  \ZZ_2\mbox{ or }1.
    \]

    Let $G_2$ and $G_3$ be a Sylow $2$-subgroup and a Sylow $3$-subgroup of $G_{\{2,3\}}$, respectively. 
    Clearly, $G_3\leqslant \bfC_{G}(G_{\{2,3\}'})$.
    Hence $G_{2'}=G_3\times G_{\{2,3\}'}=G_3\times F_{\{2,3\}'}$ is cyclic.
    Assume that $G_{2'}=\langle g\rangle=\langle g_3\rangle\times \langle g_{3'}\rangle$ such that $g_3\in G_{\{2,3\}}$.
    Then there exists $u,v\in G_2$ such that 
    \[G_2=\langle u\rangle{:}\langle v\rangle\cong \D_8\mbox{ and }G_{\{2,3\}}=\langle u^2,uv\rangle{:}(\langle g_3\rangle{:}\langle v\rangle)\cong\ZZ_2^2{:}\D_{2\times 3^e}.\]
    Note that $G=G_{\{2,3\}'}{:}G_{\{2,3\}}=\langle g_{3'}\rangle{:}G_{\{2,3\}}$ and $(G_{\{2,3\}})'=\langle u^2,uv,g_3\rangle$ has index $2$ in $G_{\{2,3\}}$.
    Then $(G_{\{2,3\}})'\leqslant \bfC_{G}(G_{\{2,3\}'})$, and hence
    \[G=\langle g_{3'}\rangle{:}G_{\{2,3\}}=(\langle u^2,uv\rangle\times \langle g_{3'}\rangle){:}(\langle g_3\rangle{:}\langle v\rangle).\]
    Note that each Sylow $r$-subgroup $R$ of $\langle g_{3'}\rangle$ is a Sylow $r$-subgroup of $G$.
    If $[R,v]=1$, then $[R,G]=1$, and so $G=R\times G_{r'}$.
    With Lemma~\ref{lem:G/N}, we deduce that $R=1$.
    Thus $\langle g_{3'}\rangle=\langle [g_{3'},v]\rangle$ and $g_{3'}^v=g_{3'}^{-1}$.
    This implies that $(g_3g_{3'})^v=(g_3g_{3'})^{-1}$. 
    Recall Lemma~\ref{lem:fit}\,(C) that $G_{\{2,3\}}\cong \ZZ_2^2{:}\D_{2\times 3^e}\cong \ZZ_{3^{e-1}}.\S_4$.
    Then $G/\langle g_{3'}g_3^3\rangle\cong \ZZ_2^2{:}\S_3\cong \S_4$ as in part~(4).
    Therefore, the proof is complete.
\end{proof}

We further determine which groups in Theorem~\ref{thm:Aut(M)} admit arc-regular triples corresponding to vertex-reversing maps satisfying the coprime condition in Problem~\ref{prob:1}.
The following lemma gives a reduction for the non-dihedral case.
Note that in the first three cases of Theorem~\ref{thm:Aut(M)}, the Hall $2'$-subgroup $G_{2'}$ is normal in $G$; that is, $G$ is $2$-nilpotent.

\begin{lemma}\label{lem:copri-G=<invos>}
    Assume that $G$ is not dihedral and $G$ can act vertex-reversingly on a map $\calM=(V,E,F)$ with $\gcd(\chi(\calM),|E|)=1$.
    Then $\calM$ is a $G$-reversing map.
    If further $G$ is $2$-nilpotent, then $G\cong \ZZ_m{:}\ZZ_2^2$ for some odd integer $m$.
\end{lemma}
\begin{proof}
    By Theorem~\ref{thm:Aut(M)}, $G_{2'}=\l g\r$ is cyclic and $G_{2}=\l u\r{:}\l v\r\cong\D_{2^{e}}$ is dihedral with $e\geqslant 2$.
    Note that $|E|=|G|/2$ is even. 
    Since $\gcd(\chi,|E|)=1$, the Euler characteristic $\chi(\calM)=|V|-|E|+|F|$ is odd.

    \textbf{Case 1.} Assume that $G$ is $2$-nilpotent, that is, $G_{2'}$ is normal.

    Let $(x,y,z)$ be a regular triple for $G$, and let $\overline{\cdot}$ be the quotient map $G\to G/G_{2'}\cong \D_{2n}$.
    If $\calM$ is $G$-bireversing, then
    \[\begin{aligned}
        \chi(\calM)\equiv \frac{|G|}{|\l x,y\r|}-\frac{|G|}{2}+\frac{|G|}{|\l x,y^z\r|}\equiv \frac{|\overline{G}|}{|\l \overline{x},\overline{y}\r|}+\frac{|\overline{G}|}{|\l \overline{x},\overline{y^z}\r|}\pmod{2}.
    \end{aligned}\]
    Since $\chi(\calM)$ is odd, $|\overline{G}:\langle\overline{x},\overline{y}\rangle|+|\overline{G}:\langle\overline{x},\overline{y^z}\rangle|$ is odd, so exactly one of the two indices is $1$.
    Hence $\overline{G}=\overline{G_2}$ equals exactly one of $\langle\overline{x},\overline{y}\rangle$ and $\langle\overline{x},\overline{y^z}\rangle$.
    This is impossible since, in a dihedral $2$-group, conjugating one generator of a generating pair yields another generating pair.
    Thus $\calM=\RevMap(G,x,y,z)$ and 
    \[\chi(\calM)\equiv \frac{|G|}{|\l x,y\r|}-\frac{|G|}{2}+\frac{|G|}{|\l x,z\r|}+\frac{|G|}{|\l y,z\r|}\equiv \frac{|\overline{G}|}{|\l \overline{x},\overline{y}\r|}+\frac{|\overline{G}|}{|\l \overline{x},\overline{z}\r|}+\frac{|\overline{G}|}{|\l \overline{y},\overline{z}\r|}\pmod{2}.\]
    This implies that 
    \begin{equation}\label{eq:rev}
        |\overline{G_2}:\langle\overline{x},\overline{y}\rangle|+
        |\overline{G_2}:\langle\overline{y},\overline{z}\rangle|+
        |\overline{G_2}:\langle\overline{x},\overline{z}\rangle|\equiv 1\pmod{2}.
    \end{equation}
    Hence at least one of the three subgroups $\l \overline{x},\overline{y}\r$, $\l \overline{x},\overline{z}\r$, and $\l \overline{y},\overline{z}\r$ equals $\overline{G_2}$.
    Without loss of generality, assume that
    \[\overline{G_2}=\langle\overline{x},\overline{y}\rangle=\langle\overline{xy}\rangle{:}\langle\overline{y}\rangle\cong\D_{2^e}.\]
    Then $\overline{z}$ is either the central involution $(\overline{xy})^{2^{e-2}}$ or a reflection $(\overline{xy})^i\overline{y}$.
    If $\overline{z}=(\overline{xy})^i\overline{y}$, then exactly one of $\langle\overline{x},\overline{z}\rangle$ and $\langle\overline{y},\overline{z}\rangle$ equals $\overline{G_2}$, making the left-hand side of~\eqref{eq:rev} even---a contradiction.
    Thus $\overline{z}=(\overline{xy})^{2^{e-2}}\in Z(\overline{G_2})$.
    If $e\geqslant 3$, then $\overline{z}\in\overline{G_2}'$, and after conjugating we may assume $z\in G_2'$.
    Since $G_{2'}$ is cyclic, $G_2/C_{G_2}(G_{2'})$ is abelian, so $[z,G_{2'}]=1$.
    As $z\in Z(G_2)$, we get $z\in Z(G)$, whence $|G:\l x,y\r|\leqslant 2$, contradicting the definition of arc-regular triples.
    Therefore $e=2$ and $G_2\cong\ZZ_2^2$.

    \textbf{Case 2.} Assume that $G$ is not $2$-nilpotent.

    We only need to prove $\calM$ is $G$-reversing.
    Suppose, to the contrary, that $\calM=\BiRevMap(G,x,y,z)$.
    Note that the edge stabilizer is isomorphic to $\ZZ_2$.
    Then each Sylow subgroup of $G$ is contained in a vertex or face stabilizer of $G$ by Lemma~\ref{lem:Sylow}\,(2).
    Hence, for $p=2$ and $3$, at least one of $\l x ,y \r$ and $\l x,y^z\r$ contains a Sylow $p$-subgroup of $G$.
    Now let $\overline{\cdot}$ be the quotient map $G\rightarrow G/\l g_{3'}g_3^3\r\cong\S_4$ as in part~(4) of Theorem~\ref{thm:Aut(M)}.
    Hence $\l \overline{x} ,\overline{y} \r$ and $\l \overline{x},\overline{y^z}\r$ are isomorphic to $\D_6$ and $\D_8$ in some ordering.
    However, a quick \textsc{Magma}~\cite{Magma} computation shows that $\S_4$ contains no such triple of involutions.
    Therefore, the proof is complete.
\end{proof}

We can now refine Theorem~\ref{thm:Aut(M)} to the groups that actually satisfy the coprime condition.

\begin{theorem}\label{thm:4Tp-coprime}
    Suppose that $\calM=(V,E,F)$ is a $G$-vertex-reversing map such that $\gcd(\chi(\calM),|E|)=1$.
    Then $G_{2'}=\l g\r$, $G_2=\l u\r{:}\l v\r$, and one of the following holds.
    \begin{enumerate} 
        \item {\rm Type I:} $G$ is a dihedral group;
        \item {\rm Type II:} $G=(\l a\r\times \l b\r){:}\l u,v\r=\l a\r{:}\l u\r \times \l b\r{:}\l v\r\cong \D_{2m}\times \D_{2n}$, where $m$ and $n$ are coprime odd integers;
        \item {\rm Type III:} $G=(\l a\r\times \l b\r \times  \l c\r){:}\l u,v\r\cong (\ZZ_m\times \ZZ_n \times \ZZ_\ell) {:}\ZZ_2^2$, where $m$, $n$, and $\ell$ are pairwise coprime odd integers, such that $(a,b,c)^u=(a^{-1},b^{-1},c)$ and $(a,b,c)^v=(a^{-1},b,c^{-1})$;
        \item {\rm Type IV:} $G=(\l u^2,uv\r\times\l g_{3'}\r){:}\l g_3,v\r$,
        where $(g_{3'}g_3)^v=(g_{3'}g_3)^{-1}$, and $G/\l g_{3'}g_3^3\r\cong\S_4$.
    \end{enumerate}
\end{theorem}
\begin{proof}
   Types I and IV are given by parts~(1) and (4) of Theorem~\ref{thm:Aut(M)}, respectively.
    The groups in parts~(2) and (3) are $2$-nilpotent, and hence $G_2\cong\ZZ_2^2$ by Lemma~\ref{lem:copri-G=<invos>}.
    If $G$ is in part~(2), then $G=(\langle g\rangle\times \langle v\rangle){:}\langle uv\rangle\cong (\ZZ_m\times \ZZ_2){:}\ZZ_2\cong\D_{4m}$, which is already covered by Type I.
    Thus, it remains to consider part~(3):
    \[
    G=\l g\r{:}\l u,v\r=(\l a\r\times\l b\r\times\l c\r){:}\langle u,v\rangle\cong (\ZZ_{m}\times\ZZ_n\times \ZZ_\ell){:}\ZZ_2^2,
    \] 
    where $bc\not=1$, $(a,b,c)^u=(a^{-1},b^{-1},c)$, and $(a,b,c)^v=(a^{-1},b,c^{-1})$.
    If exactly one of $b$ and $c$ is trivial, say $c=1$, then replacing the generating pair $\{u,v\}$ of $G_2$ by $\{uv,v\}$ yields $a^{uv}=a$, $b^{uv}=b^{-1}$, and $a^v=a^{-1}$, $b^v=b$; consequently $G\cong \l a\r{:}\l v\r \times \l b\r{:}\l uv\r\cong \D_{2m}\times\D_{2n}$, giving Type II.
    If both $b$ and $c$ are nontrivial, then $G$ is of Type III.
\end{proof}

\section{Vertex-reversing maps satisfying the coprime condition}\label{sec:final}

Let $\calM=(V,E,F)$ be a $G$-vertex-reversing map with $|V|,|F|\geqslant 3$ and Euler characteristic $\chi$, where $G$ is solvable and $\gcd(\chi,|E|)=1$.
This section is devoted to the proof of Theorem~\ref{thm:maps}.
By Theorem~\ref{thm:4Tp-coprime}, it is enough to treat separately the four possible group-theoretic types of $G$.
The corresponding classifications are given in Theorems~\ref{thm:type-I}, \ref{thm:type-II}, \ref{thm:type-III}, and~\ref{thm:type-IV}, respectively.

\subsection{Type I: dihedral groups}\

In this subsection, we determine all such $G$-vertex-reversing maps with $G=\l g\r{:}\l h\r=\D_{2n}$ in the following theorem.
The theorem is obtained by reduction from Theorem~\ref{thm:dih_map}, where the coprime condition $\gcd(\chi,|E|)=1$ is not required.

\begin{theorem}\label{thm:type-I}
    Each $G$-vertex-reversing map $\mathcal{M}=(V,E,F)$ has underlying graph $\mathbf{C}_{m}^{(n/m)}$ with $G=\langle g\rangle{:}\langle h\rangle=\D_{2n}$. 
    In addition, up to isomorphism, one of the following holds, where $\gcd(j,k,n)=1$ and $\gcd(j,n)\geq 3$.
    \vspace{-.2cm} 
    \begin{spacing}{1.2} 
    \begin{enumerate}
        \item $\calM=\RevMap(G,x,y,z)$ with $(x,y,z)=(g^{n/2},h,gh)$ and $n=2m$. 
        In this case, $\calM$ is on the projective plane with $\chi=1$ and $|F|=n/2+1$.
        \item $\calM=\RevMap(G,x,y,z)$ with $(x,y,z)=(h,g^jh,g^kh)$ and $m=\gcd(j,n)$.
        In this case, $\calM$ is always orientable and 
        \[\chi=\gcd(j,n)+\gcd(k,n)+\gcd(j-k,n)-n.\]
        \item $\calM=\BiRevMap(G,x,y,z)$ with $(x,y,z)=(h,g^jh,g^kh)$ and $m=\gcd(j,n)$. 
        In this case, $\calM$ is orientable if and only if $\gcd(j,\,2k,\,n)=2$, and 
        \[\chi=\gcd(j,n)+\gcd(2k-j,n)-n.\]
    \end{enumerate}
    \end{spacing}
\end{theorem}
\vspace{-.6cm} 

In the rest of this subsection, we assume that $\calM=(V,E,F)$ is a $G$-vertex-reversing map with $G=\l g\r{:}\l h\r =\D_{2n}$ and $(x,y,z)$ is the corresponding arc-regular triple.
First, we show that the underlying graphs of these $G$-vertex-reversing maps are cycles with multi-edges.

\begin{lemma}\label{lem:UnGraph-Dih}
    The underlying graph of $\calM$ is $\bfC_{m}^{(n/m)}$, where $m=n/|xy|$.
\end{lemma}
\begin{proof}
    Let $H=\l x,y\r$.
    Then the underlying graph is $\Gamma= \Cos(G,H,\l z\r)=(V,E)$ by Definition~\ref{def:RevMap}, and 
    \[|V|=\frac{|G|}{|H|}=\frac{2n}{2|xy|}=\frac{n}{|xy|}, \mbox{ and }|E|=\frac{|G|}{|\l z\r|}=n.\]
    Proposition~\ref{prop:CosGraph}\,(3) shows that the base graph of $\Gamma$ is $\Cos(G,H,HzH)$. 
    Note that $H\cap \l g\r$ is normal in $G$, and so it lies in the core $\Core_G(H)$ of $H$ in $G$.
    Then $|H:\Core_G(H)|\leqslant |H:H\cap \l g\r|\leqslant 2$.
    This implies that $\Cos(G,H,HzH)$ has valency $1$ or $2$.
    Thus the base graph of $\Gamma$ is either a single point or a cycle.
    Since $|V|\geqslant 3$, we have that $\Gamma\cong \bfC_{m}^{(\lambda)}$, where $m=|V|=\frac{n}{|xy|}$ and $\lambda=\frac{|E|}{|V|}=\frac{n}{m}$ as desired.
\end{proof}

The next lemma determines all arc-regular triples for dihedral groups.
\begin{lemma}\label{lem:metac-triples}
    The arc-regular triple $(x,y,z)$ is equivalent to one of the following:
    \begin{enumerate}
        \item $(g^m,h,gh)$ when $n=2m$ is even;
        \item $(g^m,h,g^2h)$ when $n=2m$ is even and $m$ is odd;
        \item $(h,g^jh,g^kh)$ with $\gcd(j,k,n)=1$ and $\gcd(j,n)\geqslant 3$.
    \end{enumerate}
\end{lemma}
\begin{proof}
    Suppose that $n=2m$ and $g^{m}\in\{x,y,z\}$.
    Note that exactly one of $x,y,z$ equals $g^m$, otherwise $G\neq\l x,y,z\r$.
    Hence we may write $\{x,y,z\}=\{g^m,g^ih,g^jh\}$ for some integers $i,j$.
    Since $\l g^m\r\lhd G$, we have $G=\l x,y,z\r=\l g^m\r \l g^ih,g^jh\r$, and then $|G:\l g^ih,g^jh\r|\leqslant |g^m|=2$.
    Then $z\neq g^m$ as $|V|=|G:\l x,y\r|\geqslant 3$.
    Recall Proposition~\ref{prop:iso-VerRevMaps} that $(x,y,z)$ is equivalent to $(y,x,z)$.
    Thus we may assume that 
    \[(x,y,z)=(g^m,g^ih,g^jh)\mbox{ for some integers $i,j$.}\]
    Note that $(g,g^ih)\mapsto (g,h)$ admits an automorphism of $G$.
    Then $(x,y,z)$ is equivalent to $(g^m,h,g^{j'}h)$ for some integer $j'$ by Proposition~\ref{prop:iso-VerRevMaps}.
    Then 
    \[G=\l x,y,z\r=\l g^m,h,g^{j'}h\r=\l g^m,g^{j'},h\r.\]
    This implies that $\gcd(m,j',n)=1$.
    Assume that $m$ is even.
    Then $\gcd(m, j',n)=\gcd(m,j',2m)=\gcd(m,j')$, and so $\gcd(n,j')=1$.
    Hence $(g,h)\mapsto (g^{j'},h)$ admits an automorphism of $G$, and so $(x,y,z)$ is equivalent to $(g^m,h,gh)$ in this case, as in Case~(1).
    Assume that $m$ is odd.
    Then $\gcd(m, j',n)=1$ implies that $\gcd(j',n)=1$ or $2$.
    If $\gcd(j',n)=1$, then  $(x,y,z)$ is equivalent to $(g^m,h,gh)$ as previous case, and so lies in Case~(1).
    If $\gcd(j',n)=2$, then $j'=2k$ and $\gcd(k,n)=1$.
    Hence $(x,y,z)$ is equivalent to $(g^m,h,g^2h)$ as $(g,h)\mapsto (g^k,h)$ admits an automorphism of $G$, and therefore Case~(2) holds.

    Now suppose that none of $x,y,z$ lies in $\l g\r$.
    As $(g,h)\mapsto (g,g^ih)$  admits an automorphism of $G$ for any integer $i$, we may assume that 
    \[(x,y,z)=(h,g^jh,g^kh)\mbox{ for some integers $j,k$.}\]
    Note that $|V|=|G:\l x,y\r|=|G:\l g^j,h\r|=|\l g\r:\l g^j\r|=\gcd(j,n)$.
    We have that $\gcd(j,n)\geqslant 3$ as $|V|\geqslant 3$.
    Since $G=\l x,y,z\r=\l h,g^jh,g^kh\r=\l g^j,g^k,h\r$, it follows that $\gcd(j,k,n)=1$, and hence Case~(3) holds.
\end{proof}

Based on the result on regular triples in the above lemma, we determine all such vertex-reversing maps.
\begin{theorem}\label{thm:dih_map}
    Each $G$-vertex-reversing map $\mathcal{M}=(V,E,F)$ has underlying graph $\mathbf{C}_{m}^{(n/m)}$ with $G=\langle g\rangle{:}\langle h\rangle=\D_{2n}$. 
    In addition, one of the following holds up to isomorphism, where $\gcd(j,k,n)=1$ and $\gcd(j,n)\geq 3$.
    \begin{center}
        \renewcommand{\arraystretch}{1.2}
        \begin{tabular}{ccccc}
        \toprule
         {\rm type} & $(x,\;y,\;z)$ &$|V|$ & $|F|$ & {\rm orientable?}\\
        \midrule
        $\Rev$ & $(g^{n/2},\;h,\;gh)$ & $n/2$ & $n/2+1$ & {\rm never}\\
        \addlinespace[2pt]
        $\Rev$ & $(g^{n/2},\;h,\;g^2 h)$ & $n/2${\rm\  is odd} & $n/2+2$ &{\rm always}\\
        \addlinespace[2pt]
        $\Rev$ & $(h,\;g^jh,\;g^kh)$ & $\gcd(j,n)$ & $\gcd(k,n)+\gcd(j-k,n)$ &{\rm always}\\
        \addlinespace[2pt]
        $\BiRev$ & $(g^{n/2},\;h,\;g h)$ & $n/2$ & $n/2$ &$n/2${\rm\  is odd}\\
        \addlinespace[2pt]
        $\BiRev$ & $(g^{n/2},\;h,\;g^2 h)$ & $n/2${\rm\  is odd} & $n/2$ &{\rm never}\\
        \addlinespace[2pt]
        $\BiRev$ & $(h,\;g^jh,\;g^kh)$ & $\gcd(j,n)$ & $\gcd(2k-j,n)$ & $\gcd(j,\,2k,\,n)=2$\\
        \bottomrule
        \end{tabular}
    \end{center}
\end{theorem}
\begin{proof}
    We consider each case in Lemma~\ref{lem:metac-triples} in turn.

    {\bf Case~1}: $(x,y,z)=(g^m,h,gh)$, where $n=2m$.

    By Lemma~\ref{lem:UnGraph-Dih}, the underlying graph is $\Cos(G,\l x,y\r,\l z\r)\cong \bfC_{m}^{(2)}$.
    If $\calM =\RevMap(G,x,y,z)$, then 
    \[|F|=\frac{|G|}{|\l x,z\r|}+\frac{|G|}{|\l y,z\r|}=\frac{n}{|xz|}+\frac{n}{|yz|}=\frac{n}{2}+1.\]
    Then $\chi(\calM)=|V|-|E|+|F|=n/2-n+n/2+1=1$, and hence the supporting surface is the projective plane and is non-orientable, as in row~1.
    If $\calM =\BiRevMap(G,x,y,z)$, then 
    \[|F|=\frac{|G|}{|\l x,y^z\r|}=\frac{n}{|xy^z|}=\frac{n}{2}.\]
    By Lemma~\ref{cri-Orient}, $\calM$ is orientable if and only if 
    \[G^+:=\l xy,z,xzy\r=\l g^{\frac{n}{2}}h,gh,g^{\frac{n}{2}+1}\r= \l g^{\frac{n}{2}+1}\r{:}\l gh\r \]
    has index $2$ in $G$.
    As $\gcd(1+n/2,n)=2$ if and only if $n/2$ is odd, we have that $\calM$ is orientable if and only if $n/2$ is odd as in row~4.
    
    {\bf Case~2}: $(x,y,z)=(g^m,h,g^2h)$, where $n=2m$ and $m$ is odd.

    Lemma~\ref{lem:UnGraph-Dih} implies that the underlying graph is $\Cos(G,\l x,y\r,\l z\r)\cong \bfC_{m}^{(2)}$.
    By an argument similar to that in the previous case, we conclude that $\RevMap(G,x,y,z)$ has $2+m$ faces and $\BiRevMap(G,x,y,z)$ has $m$ faces.
    By Lemma~\ref{cri-Orient}, $\calM$ is orientable if and only if 
    \[G^+=\begin{cases}
        \l xy,xz\r=\l g^2\r{:}\l g^mh\r, & \text{if }\calM=\RevMap(G,x,y,z); \\
        \l xy,z,xzy\r=\l g^{m+2},g^{m-2}\r{:}\l g^2h\r, & \text{if }\calM=\BiRevMap(G,x,y,z)
    \end{cases} \]
    has index $2$ in $G$.
    Since $n=2m$ and $m$ is odd, we have that $\l g^2\r<\l g\r$ and $\l g^{m+2},g^{m-2}\r=\l g\r$.
    It follows that $\RevMap(G,x,y,z)$ is orientable and $\BiRevMap(G,x,y,z)$ is non-orientable, as in rows~2 and 5.

    {\bf Case~3}: $(x,y,z)=(h,g^jh,g^kh)$, where $\gcd(j,k,n)=1$ and $\gcd(j,n)\geqslant 3$.

    Recall Lemma~\ref{lem:UnGraph-Dih} that the underlying graph is $\Cos(G,\l x,y\r,\l z\r)\cong \bfC_{m}^{(n/m)}$ with $m=n/|xy|=n/|g^j|=\gcd(n,j)$.
    Note that
    \begin{enumerate}
        \item[(a)] if $\calM=\RevMap(G,x,y,z)$, then 
        \[|F|=\frac{|G|}{|\l x,z\r|}+\frac{|G|}{|\l y,z\r|}=\frac{n}{|xz|}+\frac{n}{|yz|}=\frac{n}{|g^k|}+\frac{n}{|g^{j-k}|}=\gcd(k,n)+\gcd(j-k,n);\]
        \item[(b)] if $\calM=\BiRevMap(G,x,y,z)$, then 
        \[|F|=\frac{|G|}{|\l x,y^z\r|}=\frac{n}{|xy^z|}=\frac{n}{|g^{2k-j}|}=\gcd(2k-j,n).\]
    \end{enumerate}
    By Lemma~\ref{cri-Orient}, $\calM$ is orientable if and only if 
    \[G^+=\begin{cases}
        \l xy,xz\r=\l g^j,g^k\r, & \text{if }\calM=\RevMap(G,x,y,z); \\
        \l xy,z,xzy\r=\l g^{j},g^{2k-j}\r {:}\l g^k h\r, & \text{if }\calM=\BiRevMap(G,x,y,z)
    \end{cases} \]
    has index $2$ in $G$.
    Thus $\RevMap(G,x,y,z)$ is orientable, and $\BiRevMap(G,x,y,z)$ is orientable if and only if $\gcd(j,2k-j,n)=\gcd(j,2k,n)=2$ as in rows~3 and 6.
\end{proof}

Using the above theorem, we now obtain the classification stated in Theorem~\ref{thm:type-I}.
\begin{proof}[{\bf Proof of Theorem~\ref{thm:type-I}}]
It is straightforward to see that the Euler characteristics of the maps in rows~2, 4 and 5 in the table of Theorem~\ref{thm:dih_map} are $2$, $0$ and $0$, respectively.
For rows~2 and 4, $|E|=n$ is even, and hence $\chi=2$ is not coprime to $|E|$.
For row~5, $\chi=0$ is also not coprime to $|E|$.
Consequently, $\mathcal{M}$ can only consist of maps in rows~1, 3 and 6, which are exactly the three cases in Theorem~\ref{thm:type-I}.

It remains to verify that each of the three remaining rows can indeed be realized under the coprime condition and the assumptions $|V|,|F|\geqslant 3$.

First, consider row~1.
For any even integer $n\geqslant 6$, we have $|V|=n/2\geqslant 3$, $|F|=n/2+1\geqslant 4$, and $\chi=1$ is coprime to $|E|=n$.

For row~3, take $n=pq$, $j=p$, and $k=q$, where $p<q$ are odd primes with $q\not\equiv -1\pmod p$.
Then $\gcd(j,k,n)=1$, $|V|=\gcd(j,n)=p\geqslant 3$, and $|F|=\gcd(k,n)+\gcd(j-k,n)=q+1\geqslant 3$.
Moreover, $\chi=p+q+1-pq$, which is coprime to $|E|=pq$ by the choice of $p$ and $q$.

For row~6, take $n=pq$, $j=p$, and $k=(p+q)/2$, where $p$ and $q$ are distinct odd primes.
Then $\gcd(j,k,n)=1$, $|V|=\gcd(j,n)=p\geqslant 3$, and $|F|=\gcd(2k-j,n)=\gcd(q,pq)=q\geqslant 3$.
Moreover, $\chi=p+q-pq$, which is coprime to $|E|=pq$.
Hence all three rows can indeed be realized under the coprime condition.
\end{proof}

\subsection{Type II}\

In this subsection, we determine all $G$-vertex-reversing maps $\calM=(V,E,F)$ with $|V|,|F|\geqslant 3$ and $\gcd(\chi(\calM),|E|)=1$ such that $G$ is of Type~II in Theorem~\ref{thm:4Tp-coprime}:
\[G=\l a\r{:}\l u\r\times \l b\r{:}\l v\r\cong \D_{2m}\times\D_{2n},\]
where $m,n>1$ are coprime odd integers.

\begin{theorem}\label{thm:type-II}
    Let $w=uv$.
    Then $\calM$ is non-orientable and $\chi(\calM)=m+n-mn$.
    In addition, up to isomorphism, $\calM$ is one of the following:
    \begin{enumerate}
    	\item $\RevMap(G,u,v,abw)$ with underlying graph $\bfC_m\times\bfC_n$;
    	\item $\RevMap(G,u,abw,v)$ with underlying graph $\bfC_n^{(2m)}$;
    	\item $\RevMap(G,v,abw,u)$ with underlying graph $\bfC_m^{(2n)}$.
    \end{enumerate}
    Moreover, each of the above three maps satisfies $\gcd(\chi(\calM),|E|)=1$.
\end{theorem}
\begin{proof}
    Let $\chi=\chi(\calM)$.
    By Lemma~\ref{lem:copri-G=<invos>}, $\calM=\RevMap(G,x,y,z)$ for some arc-regular triple $(x,y,z)$, and then
    \[\chi=\frac{|G|}{2|xy|}+\frac{|G|}{2|xz|}+\frac{|G|}{2|yz|}-\frac{|G|}{2}=\frac{2mn}{|xy|}+\frac{2mn}{|xz|}+\frac{2mn}{|yz|}-2mn.\]
    
    First, we show that $(x,y,z)$ is equivalent to some ordering of $\{u,v,abw\}$.
    Let $\overline{\cdot}$ be the natural quotient map $G\to G/\l a,b\r$.
    Then $\overline{G}=\l\overline{x},\overline{y},\overline{z}\r\cong\ZZ_2\times\ZZ_2$ and $\overline{x},\overline{y},\overline{z}\neq 1$ as $|\l a,b\r|=mn$ is odd.
    If $|\{\overline{x},\overline{y},\overline{z}\}|=2$, then the three values $|\overline{xy}|, |\overline{xz}|, |\overline{yz}|$ are $2$, $2$, and $1$ in some ordering.
    Since $|G|_2=4$, we have
    \[\chi\equiv \frac{2}{|\l xy\r_2|}+\frac{2}{|\l xz\r_2|}+\frac{2}{|\l yz\r_2|}\equiv 1+1+2\equiv 0\pmod{2}.\]
    Then $2\mid \gcd(\chi,|E|)$, a contradiction.
    Thus $|\{\overline{x},\overline{y},\overline{z}\}|=3$, and so 
    \[\{x,y,z\}=\{a^{i_1}b^{i_2}u,\;a^{j_1}b^{j_2}v,\;a^{k_1}b^{k_2}w\},\mbox{ for some integers $i_1,i_2,j_1,j_2,k_1,k_2.$}\]
    As $|x|=|y|=|z|=2$, we have that $i_2=j_1=0$.
    Note that $(a,b,u,v)\mapsto (a,b,a^{i_1}u,b^{j_2}v)$ admits an automorphism of $G$.
    Then we may assume that 
    \[\{x,y,z\}=\{u,\;v,\;a^{k_1}b^{k_2}w\},\mbox{ for some integers $k_1,k_2.$}\]
    As $G=\l x,y,z\r=\l u,v,a^{k_1}b^{k_2}w\r=\l a^{k_1},b^{k_2},u,v\r$, we have that $\gcd(k_1,m)=\gcd(k_2,n)=1$.
    Then $(a,b,u,v)\mapsto (a^{k_1},b^{k_2},u,v)$ also admits an automorphism of $G$, and so $(x,y,z)$ is equivalent to some ordering of $\{u,v,abw\}$.
    Then 
    \[\chi=\frac{2mn}{|xy|}+\frac{2mn}{|xz|}+\frac{2mn}{|yz|}-2mn=m+n-mn.\]
    Since $\chi$ is odd, $\calM$ is non-orientable.
    Then 
    \[\gcd(\chi,|E|)=\gcd(m+n-mn,2mn)=\gcd(m+n,mn).\]
    Since $\gcd(m,n)=1$, for each prime divisor $p$ of $mn$, $p$ divides exactly one of $m$ and $n$.
    Then $\gcd(m+n,mn)=1$, and hence $\gcd(\chi,|E|)=1$.
    Thus it remains only to determine the underlying graph $\Gamma=\Cos(G,H,J)$ of $\calM$, where $H=\l x,y\r$ and $J=\l z\r$.
    Proposition~\ref{prop:iso-VerRevMaps} shows that $(x,y,z)$ is equivalent to $(y,x,z)$.
    Then $(x,y,z)$ is equivalent to either 
    \[(u,v,abw),\ (u,abw,v)\mbox{ or }(v,abw,u).\]

    \textbf{Case 1.} Assume that $(x,y,z)=(u,v,abw)$.

    Let $K=H\cap H^z$.
    Then $K=\l u,v\r\cap \l a^2u, b^2v\r$ is trivial as $n,m>1$.
    By Proposition~\ref{prop:CosGraph}, we have that $\Gamma$ is simple and $\Gamma\cong \Cos(G,H,HzH)$.
    Recall that $G=\l a\r{:}\l u\r\times \l b\r{:}\l v\r\cong \D_{2m}\times\D_{2n}$ and $H=\l u\r \times \l v\r$.
    By~Lemma~\ref{lem:product-coset-graph}, $\Gamma$ is isomorphic to $\bfC_m\times \bfC_n$.
    
    \textbf{Case 2.} Assume that $(x,y,z)=(u,abw,v)$.

    Clearly, $|V|=|G:\l x,y\r|=n$ and $|E|=|G|/2=2mn$.
    By Proposition~\ref{prop:CosGraph}, the edge multiplicity is $2m$ and the edge-valency is $|H:H\cap J|=4m$; hence the vertex-valency is $2$.
    Thus $\Gamma\cong\bfC_n^{(2m)}$.
    
    \textbf{Case 3.} Assume that $(x,y,z)=(v,abw,u)$.

    Repeating the process in the previous case, we obtain that $\Gamma\cong\bfC_m^{(2n)}$.
    The proof is now complete by exhaustion of all possibilities.
\end{proof}

\subsection{Type III}\

In this subsection, we determine all the $G$-vertex-reversing maps $\calM=(V,E,F)$ with $|V|,|F|\geqslant 3$ and $\gcd(\chi(\calM),|E|)=1$ such that $G$ is of Type~III in Theorem~\ref{thm:4Tp-coprime}:
\[
	G=(\l a\r\times\l b\r\times \l c\r){:}(\l u,v\r)\cong (\ZZ_m\times \ZZ_n\times \ZZ_\ell) {:} \D_4,
\]
where $m,n,\ell>1$ are pairwise coprime odd integers and 
\[
	(a,b,c)^u=(a^{-1},b^{-1},c),
	\qquad
	(a,b,c)^v=(a^{-1},b,c^{-1}).
\]

\begin{theorem}\label{thm:type-III}
    Let $w=uv$.
    Then $\calM$ is non-orientable and $\chi(\calM)=mn+m\ell+n\ell-2mn\ell$.
    In addition, up to isomorphism, $\calM$ is one of the following:
    \begin{enumerate}
    	\item $\RevMap(G,u,av,bcw)$ with underlying graph $(\bfC_n\times\bfC_\ell)^{(m)}$;
    	\item $\RevMap(G,u,bcw,av)$ with underlying graph $(\bfC_m\times\bfC_\ell)^{(n)}$;
    	\item $\RevMap(G,av,bcw,u)$ with underlying graph $(\bfC_m\times\bfC_n)^{(\ell)}$.
    \end{enumerate}
    Moreover, each of the above three maps satisfies $\gcd(\chi(\calM),|E|)=1$.
\end{theorem}

\begin{proof}
    Lemma~\ref{lem:copri-G=<invos>} implies that $\calM=\RevMap(G,x,y,z)$ for some arc-regular triple $(x,y,z)$.
    Let $w=uv$, and let $\overline{\cdot}$ be the natural quotient map $G\to G/\l a,b,c\r\cong\ZZ_2^2$.
    Then $\overline{G}=\l\overline{x},\overline{y},\overline{z}\r\cong\ZZ_2^2$.
    By arguments similar to those in the proof of Theorem~\ref{thm:type-II}, we conclude that $\{\overline{x},\overline{y},\overline{z}\}=\{\overline{u},\overline{v},\overline{w}\}$ and
    \[
    \{x,y,z\}=\{a^{i_1}b^{i_2} u,\;a^{j_1}c^{j_2} v,\;b^{k_1}c^{k_2} w\}\mbox{ for some integers $i_1,i_2,j_1,j_2,k_1,k_2$.}
    \]
    Note that $(a,b,c,u,v)\mapsto (a,b,c,a^{i_1}b^{i_2} u, a^{i_1}c^{j_2} v)$ admits an automorphism of $G$.
    We may assume that $i_1=i_2=j_2=0$, and hence
    \[\{x,y,z\}=\{u,\;a^{j_1}v,\;b^{k_1}c^{k_2} w\}\mbox{ for some integers $j_1,k_1,k_2$.}\]
    Recall that $G=\l x,y,z\r$.
    Then 
    \[\begin{aligned}
        \l a,b,c\r{:}\l u,v\r&=\l x,y,z\r=\l u,\; a^{j_1}v,\;b^{k_1}c^{k_2} w \r=\l u,\; a^{j_1}v,\; a^{j_1} w,\; b^{k_1}c^{k_2} w\r\\
        &= \l a^{j_1},\;u,\; v,\; w,\; b^{k_1}c^{k_2} w\r=\l a^{j_1},b^{k_1},c^{k_2}\r{:}\l u,v\r.
    \end{aligned}\]
    This implies that $\gcd(j_1,m)=\gcd(k_1,n)=\gcd(k_2,\ell)=1$.
    Then $(a,b,c,u,v)\mapsto (a^{j_1},b^{k_1},c^{k_2},u,v)$ admits an automorphism of $G$.
    We have that $(x,y,z)$ is equivalent to some ordering of $\{ u,av,bcw\}$.
    Let $\chi=\chi(\calM)$.
    Then 
    \[\chi=\frac{2mn\ell}{|xy|}+\frac{2mn\ell}{|xz|}+\frac{2mn\ell}{|yz|}-2mn\ell=mn+m\ell+n\ell-2mn\ell.\]
    Clearly, $\chi$ is odd, and hence $\calM$ is non-orientable.
    Then 
    \[\gcd(\chi,|E|)=\gcd(mn+m\ell+n\ell-2mn\ell,4mn\ell )=\gcd(mn+m\ell+n\ell,mn\ell).\]
    Since $m$, $n$, and $\ell$ are pairwise coprime odd integers, for each prime divisor $p$ of $mn\ell$, $p$ divides exactly one of $m$, $n$ and $\ell$.
    Then $\gcd(mn+m\ell+n\ell,mn\ell)=1$, and hence $\gcd(\chi,|E|)=1$.
    Thus it remains only to determine the underlying graph $\Gamma=\Cos(G,H,J)$ of $\calM$, where $H=\l x,y\r$ and $J=\l z\r$.
    Proposition~\ref{prop:iso-VerRevMaps} shows that $(x,y,z)$ is equivalent to $(y,x,z)$.
    Then $(x,y,z)$ is equivalent to either 
    \[(u,av,bcw),\ (u,bcw,av)\mbox{ or }(av,bcw,u).\]
    Since the arguments for the three cases are essentially the same, we only present the first one: assume that $(x,y,z)=(u,av,bcw)$.
    
    Let $H=\l x,y\r=\l u,av\r$, $J=\l z\r=\l bcw\r$, and let $K=H\cap H^z=\l a\r$.
    By Proposition~\ref{prop:CosGraph}, the underlying graph of $\calM$ is isomorphic to $\Gamma=\Cos(G,H,J)$ with edge-multiplicity $|K:H\cap J|=|K|=m$, and the base graph $\Gamma_0$ of $\Gamma$ is isomorphic to $\Cos(G,H,HzH)$.
    Note that $\Core_G(H)=\l a\r$.
    Then $\Gamma_0\cong \Cos(G_0,H_0,H_0zH_0)$, where 
    \[G_0=\l b\r{:}\l u\r \times \l c\r{:}\l v\r\cong\D_{2n}\times\D_{2\ell}\mbox{ and }H_0=\l u\r\times\l v\r.\]
    By Lemma~\ref{lem:product-coset-graph}, $\Gamma_0$ is isomorphic to $\bfC_n\times \bfC_\ell$.
\end{proof}

\subsection{Type IV and the proof of Theorem~\ref{thm:maps}}\

In this subsection, we determine all $G$-vertex-reversing maps  with $|V|,|F|\geqslant 3$ and $\gcd(\chi(\calM),|E|)=1$ such that $G$ is of Type~IV in Theorem~\ref{thm:4Tp-coprime}:
\[
    G=\l u^2,uv\r{:}(\l g\r{:}\l v\r)=(\l u^2,uv\r\times\l g_{3'}\r){:}\l g_3,v\r\cong (\ZZ_2^2\times\ZZ_n){:}\D_{2\cdot 3^{f+1}}\cong \ZZ_{3^fn}.\S_4,
\]
where $f\geqslant 0$, $\gcd(6,n)=1$, $G_2=\l u\r{:}\l v\r\cong\D_8$ is a Sylow $2$-subgroup of $G$, $G_{2'}=\l g\r$ is a Hall $2'$-subgroup of $G$, and $g^v=g^{-1}$.
Note that $G/\l g^3\r\cong\S_4$, and then $g$ rotates the three involutions of $\langle uv, u^2\rangle\cong\ZZ_2^2$.
We may assume that
\[
(uv,\; u^3v,\; u^2)^g = (u^3v,\; u^2,\; uv).
\]

\begin{lemma}\label{lem:Cn.S4-coprime-triples}
Assume that $\calM$ is a $G$-vertex-reversing map with $|V|,|F|\geqslant 3$ and $\gcd(\chi(\calM),|E|)=1$.  Then $\calM\cong \RevMap(G,x,y,z)$, where 
\[
    \{x,y,z\}=\{v,\ g^{1+3i}v,\ u^2\}\mbox{ with $\l g^{1+3i}\r=\l g\r$.}
\]
Moreover, $\chi(\calM)=4-3^{f+1}\cdot n$ and $|E|=4\cdot 3^{f+1}\cdot n$.
\end{lemma}
\begin{proof}
    By Lemma~\ref{lem:copri-G=<invos>}, $\calM\cong \RevMap(G,x,y,z)$ for some arc-regular triple $(x,y,z)$ of $G$.
    Clearly, $|E|=|G|/2=4|g|=4\cdot 3^{f+1}\cdot n$ is even.

    Let $\overline{\cdot}$ be the natural quotient map $G\rightarrow G/\l g^3\r\cong \S_4$.
    Since $|g^3|$ is odd, we have $\overline{x}$, $\overline{y}$, and $\overline{z}$ are involutions in $\overline{G}$; and
    \[\begin{aligned}
        \chi(\calM)&\equiv |V|+|F|-|E|\equiv\frac{|G|}{|\l x,y\r|}+\frac{|G|}{|\l x,z\r|}+\frac{|G|}{|\l y,z\r|}&\pmod{2}\\
        &\equiv\frac{|\overline{G}|}{|\overline{\l x,y\r}|}+\frac{|\overline{G}|}{|\overline{\l x,z\r}|}+\frac{|\overline{G}|}{|\overline{\l y,z\r}|}\equiv \frac{12}{|\overline{x}\overline{y}|}+\frac{12}{|\overline{x}\overline{z}|}+\frac{12}{|\overline{y}\overline{z}|}&\pmod{2}.
    \end{aligned}\]
    Since $\gcd(\chi(\calM),|E|)=1$, $\chi(\calM)$ is odd, and hence $12/|\overline{x}\overline{y}|+12/|\overline{x}\overline{z}|+12/|\overline{y}\overline{z}|$ is odd.
    Note that \(\langle\overline{x},\overline{y},\overline{z}\rangle=\overline{G}\cong\S_4\).
    A direct computation with \textsc{Magma}~\cite{Magma} shows that there are exactly $24$ possible such $3$-element sets $\{\overline{x},\overline{y},\overline{z}\}$ in $\overline{G}$, and for each of them there exists an element $\gamma\in G$ such that
    \[
    \{\overline{v},\overline{gv},\overline{u}^2\}=\{\overline{x},\overline{y},\overline{z}\}^{\overline{\gamma}}
    = \{\overline{x^\gamma},\overline{y^\gamma},\overline{z^\gamma}\}.
    \]
    By Proposition~\ref{prop:iso-VerRevMaps}, we may assume that $\{\overline{x},\overline{y},\overline{z}\}=\{\overline{v},\overline{gv},\overline{u}^2\}$.
    Then 
    \[\{x,\;y,\;z\}=\{g^{3i}v,\;g^{1+3j}v,\;u^2\}\mbox{ for some integers $i,j$}.\]
    Since $|g|$ is odd, there exists $h\in\l g\r$ such that $h^2=g^{3i}$.
    Then 
    \[\{x,\;y,\;z\}^h=\{g^{3i}v,\;g^{1+3j}v,\;u^2\}^h=\{v,\; g^{1+3(j-i)}v,\; u^2\}.\]
    Hence $(x,y,z)$ is equivalent to some ordering of $\{v,g^{1+3(j-i)}v,u^2\}$.
    Then 
    \[\begin{aligned}
        G&=\l x,y,z\r=\l v,g^{1+3(j-i)}v,u^2\r=\l g^{1+3(j-i)},v,u^2\r=\l g^{1+3(j-i)},v,u^2,(u^2)^{g^{1+3(j-i)}v}\r\\
        &=\l g^{1+3(j-i)}v,v,u^2,uv\r=\l g^{1+3(j-i)}v,u ,v \r.
    \end{aligned}\]
    This implies that $\l g^{1+3(j-i)}\r=\l g\r$.
    Thus we conclude that $(x,y,z)$ is equivalent to some ordering of $\{v,g^{1+3i}v,u^2\}$ for some integer $i$ such that $\l g\r=\l g^{1+3i}\r$.

    Now we calculate $\chi(\calM)$.
    Note that 
    \[\begin{aligned}
        \chi(\calM)&= |V|+|F|-|E|=\frac{|G|}{|\l x,y\r|}+\frac{|G|}{|\l x,z\r|}+\frac{|G|}{|\l y,z\r|}-\frac{|G|}{2}\\
        &=|G|\left(\frac{1}{2|xy|}+\frac{1}{2|xz|}+\frac{1}{2|yz|}-\frac{1}{2}\right)\\
        &=4\cdot 3^{f+1}\cdot n\cdot \left(\frac{1}{|g^{1+3i}|}+\frac{1}{|vu^2|}+\frac{1}{|g^{1+3i}vu^2|}-1\right)
    \end{aligned}\]
    Clearly, $|g^{1+3i}|=|g|=3^{f+1}\cdot n$ and $|vu^2|=2$.
    Then $|gvu^2|$ is even, and so 
    \[
        \begin{aligned}
            |gvu^2|&=2|gvu^2\cdot gvu^2|=2|gvu^2vg^{-1}u^2|=2|gvu^2vuvg^{-1}|\\
            &=2|vu^2vuv|=2|uv|=4.
        \end{aligned}
    \]
    Since $gvu^2$ maps $g^3$ to $g^{-3}$ by conjugation, we have that
    \begin{equation}\label{eq:xy}
        |g^{1+3i}vu^2|=|gvu^2|=4
    \end{equation}
    Therefore, $\chi(\calM)=4-3^{f+1}\cdot n$.
\end{proof}

Now we are ready to determine all $G$-vertex-reversing maps in this type.

\begin{theorem}\label{thm:type-IV}
    $\calM$ is non-orientable with $\chi(\calM)=4-3^{f+1}\cdot n$.
    In addition, up to isomorphism, $\calM$ is one of the following, where $\l g^{1+3i}\r=\l g\r$.
    \begin{enumerate}
        \item $\RevMap(G,\;v,\;g^{1+3i}v,\;u^2)$, with underlying graph $\K_4^{(2\cdot 3^{f}\cdot n)}$;
        \item $\RevMap(G,\;v,\;u^2,\;g^{1+3i}v)$, with underlying graph $\bfC_{3^{f+1}\cdot n}[\overline{\K_2}]$;
        \item $\RevMap(G,\;g^{1+3i}v,\;u^2,\;v)$, with underlying graph $\bfC_{3^{f+1}\cdot n}^{(4)}$.
    \end{enumerate}
    Moreover, each of the above three maps satisfies $\gcd(\chi(\calM),|E|)=1$.
\end{theorem}
\begin{proof}
Lemma~\ref{lem:Cn.S4-coprime-triples} shows that $\chi(\calM)=4-3^{f+1}\cdot n$, $|E|=4\cdot 3^{f+1}\cdot n$ and $\calM\cong \RevMap(G,x,y,z)$ for some ordering $(x,y,z)$ of $\{v, g^{1+3i} v, u^2\}$, where $\l g\r =\l g^{1+3i}\r$.
Hence, for each ordering, we always have that $\chi(\calM)$ is odd and
\[\begin{aligned}
    \gcd(\chi(\calM),|E|)&=\gcd(4-3^{f+1}\cdot n,4\cdot 3^{f+1}\cdot n)=\gcd(4-3^{f+1}\cdot n,3^{f+1}\cdot n)\\
    &=\gcd(4,3^{f+1}\cdot n)=1.
\end{aligned}\]
By Proposition~\ref{prop:iso-VerRevMaps}, it remains only to determine the underlying graph $\Gamma$ of $\mathcal{M}=\RevMap(G,x,y,z)$ with $(x,y,z)$ running over
\[(v,\;g^{1+3i}v,\;u^2),\quad (v,\;u^2,\;g^{1+3i}v),\mbox{ and }(g^{1+3i}v,\;u^2,\;v).\]
Note that $\Gamma=\Cos(G,H,\l z\r)$ has base graph $\Gamma_0=\Cos(G,H,HzH)$ with edge-multiplicity $|H\cap H^z:H\cap \l z\r|$, where $H=\l x,y\r$ by Proposition~\ref{prop:CosGraph}.

\textbf{Case 1.} Assume that $(x,y,z)=(v,\;g^{1+3i}v,\;u^2)$.

In this case, $H=\l v,g^{1+3i}v\r=\l g\r {:}\l v\r$.
Then $\l g^3\r \lhd H$ is normal in $G$.
Since $G/\l g^3\r\cong\S_4$ and $H/\l g^3\r\cong \S_3$, we deduce that $\Gamma_0=\K_4$.
Then the edge-multiplicity of $\Gamma$ is $|E|/6=2\cdot 3^{f}\cdot n$, and so $\Gamma\cong \K_4^{(2\cdot 3^{f}\cdot n)}$ as in part~(1).

\textbf{Case 2.} Assume that $(x,y,z)=(v,\;u^2,\;g^{1+3i}v)$.

In this case, $H=\l v,u^2\r\cong \ZZ_2^2$.
Note that $z$ does not normalize $H$ and $\bfC_{H}(z)=1$.
It follows that $H\cap H^z=1$, and hence $\Gamma=\Gamma_0$ is a simple graph of valency $4$.
Let $N=\langle u^2,uv\rangle\lhd G$.
Each orbit of $N$ on $V$ has length $|N:N\cap H|= 2$. 
It follows that $N$ has $|V|/2=|G:H|/2=3^{f+1}\cdot n$ orbits on $V$.
Define $\Gamma/N$ to be the simple graph whose vertices are the $N$-orbits on $V\Gamma$, and two orbits are adjacent whenever they have adjacent vertices in $\Gamma$.
Then $\Gamma/N\cong \Cos(\overline{G},\overline{H},\overline{HzH})$ where $\overline{\cdot}$ is the natural quotient map $G\rightarrow G/N$.
Note that $\overline{H}\cong \ZZ_2$.
Then $\Gamma/N$ has valency $2$, and hence $\Gamma/N\cong \bfC_{3^{f+1}\cdot n}$.
Recall that each orbit of $N$ has length $2$ and $\Gamma$ has valency $4$.
This implies that $\Gamma\cong \bfC_{3^{f+1}\cdot n}[\overline{\K_2}]$ as in part~(2).

\textbf{Case 3.} Assume that $(x,y,z)=(g^{1+3i}v,\;u^2,\;v)$.

Then by Equation~\ref{eq:xy}, we have $|H|=2|xy|=8$, and so $H\cong\D_8$. 
Note that $y^x=(uv)^v=u^3v$, and hence $\l uv, u^2\r\leqslant H$.
Thus either $H\cap H^z=\l uv, u^2\r\cong \ZZ_2^2$ or $z$ normalizes $H$.
The latter case is impossible since $|V|=|G:H|=3^{f+1}\cdot n>2$.
Hence $|V|=|G:H|=3^{f+1}\cdot n$ and the edge-multiplicity equals $|H\cap H^z: H\cap \l z\r|=|H\cap H^z|=4$.
Then $\Gamma_0$ has $|G|/8=3^{f+1}\cdot n$ edges.
This implies that $\Gamma_0=\bfC_{3^{f+1}\cdot n}$.
Hence $\Gamma=\bfC_{3^{f+1}\cdot n}^{(4)}$ as in part~(3).
\end{proof} 

The preceding four subsections have established the classification for each of the four types.
To complete the proof of Theorem~\ref{thm:maps}, we now combine these results.

\begin{proof}[{\bf Proof of Theorem~\ref{thm:maps}}]
    By Lemma~\ref{lem:copri-G=<invos>}, either $\calM$ is $G$-bireversing (and $G$ is dihedral), or $\calM$ is $G$-reversing.  
    For $G$-bireversing maps, Theorem~\ref{thm:type-I}\,(3) gives exactly the map described in Case~1 of the theorem.
    For $G$-reversing maps, Theorem~\ref{thm:4Tp-coprime} shows that $G$ must belong to one of the four group-theoretic types I-IV.
    Applying parts~(1) and (2) of Theorem~\ref{thm:type-I} and Theorems~\ref{thm:type-II},~\ref{thm:type-III},~\ref{thm:type-IV} in each case, we obtain exactly the entries in Table~\ref{tab:classification}.
    Thus the two cases cover all possibilities, and the proof is complete.
\end{proof}

\vskip0.1in
\noindent\thanks{{\bf Acknowledgments.}
We are also grateful to the anonymous referees for their valuable comments and suggestions, which have greatly improved the paper.
This work was partially supported by NSFC: 12350710787.}

\end{document}